\documentclass[titlepage,11pt]{article}

\usepackage[dvips]{graphicx}

\usepackage[myheadings]{fullpage}
\usepackage{pmetrika}
\usepackage{pmbib}

\usepackage{submit}

\usepackage[english]{babel}
\usepackage[utf8]{inputenc}
\usepackage[T1]{fontenc}
\usepackage{amsmath,amssymb,amsfonts}
\usepackage{bbm}
\usepackage{nicefrac}
\usepackage{tikz}
\usepackage{subfigure}
\usepackage{dashrule}
\newcommand{\bs}[1]{\boldsymbol{#1}}
\def\Erw{\mathbb{E}}
\def\diag{\mathrm{diag}}
\def\vec{\mathrm{vec}}
\def\Cov{\mathrm{cov}}
\def\tr{\mathrm{tr}}

\def\var{\mathrm{var}}
\def\Rel{\mathrm{Rel}}
\def\bs{\boldsymbol}

\def\bX{\bs{X}}
\def\bep{\bs{\varepsilon}}
\def\bmX{\bs{\overline{X}}}
\def\bnull{\bs{0}}
\def\bSigma{\bs{\Sigma}}
\def\Ex{\mathop{\mathrm{I\!E}}\nolimits}
\newcommand{\R}{\mathbb{R}}

\newcommand{\dist}{\mbox{$\, \stackrel{d}{\longrightarrow} \,$}}
\newcommand{\prob}{\mbox{$\, \stackrel{p}{\longrightarrow} \,$}}

\newcommand{\kone}{{\mathbf{1}}_k}

\usepackage{color}
\definecolor{dkgreen}{rgb}{0,0.6,0}
\definecolor{gray}{rgb}{0.5,0.5,0.5}
\definecolor{mauve}{rgb}{0.58,0,0.82}
 
\usepackage{listings}
\begin{document}

\begin{titlepage}

\title{Resampling-based Inference Methods for Comparing Two Coefficient Alpha}


\author{Markus Pauly and Maria Umlauft}

\affil{Institute of Statistics, Ulm University, Germany}

\author{Ali Ünlü}
\affil{Technical University of Munich, Germany}

\vspace{\fill}\centerline{\today}\vspace{\fill}

\thanks{The work of Markus Pauly and Maria Umlauft was supported by the German Research Foundation project DFG-PA 2409/3-1.}
\linespacing{1}
\contact{Correspondence should be sent to\\

\noindent E-Mail: markus.pauly@uni-ulm.de \\
\noindent Phone: (+49) 731/ 50-33105 \\
\noindent Fax: (+49) 731/50-33110 \\
\noindent Website: uni-ulm.de/mawi/statistics/team }

\end{titlepage}

\setcounter{page}{2}
\vspace*{2\baselineskip}

\RepeatTitle{Resampling-based Inference Methods for Comparing Two Coefficient Alpha}\vskip3pt

\linespacing{2.0}

\abstracthead
\begin{abstract}
The two-sample problem for Cronbach's coefficient $\alpha_C$, as an estimate of test or composite score reliability, has attracted little attention, 
compared to the extensive treatment of the one-sample case. It is necessary to compare the reliability of a test for different subgroups, for different tests or the short and long 
 forms of a test. In this paper, we study statistically how to compare two coefficients $\alpha_{C,1}$ and $\alpha_{C,2}$.
The null hypothesis of interest is $H_0 : \alpha_{C,1} = \alpha_{C,2}$, which we test against one-or two-sided alternatives. 
For this purpose, resampling-based permutation and bootstrap tests are proposed. 
These statistical tests ensure a better control of the type I error, in finite or very small sample sizes, when the state-of-affairs \textit{asymptotically distribution-free} (ADF) 
large-sample test may fail 
to properly attain the nominal significance level.
We introduce the permutation and bootstrap tests for the two-group multivariate non-normal models under the general ADF setting, thereby improving on the small sample properties of the well-known ADF asymptotic test.
By proper choice of a studentized test statistic, the resampling tests are modified such that they are still asymptotically valid, 
if the data may not be exchangeable. 
The usefulness of the proposed resampling-based testing strategies is demonstrated in an extensive simulation study and illustrated 
by real data applications.

\end{abstract}\vspace{\fill}\pagebreak

\section{Introduction}\label{sec:intro}
Reliability is a cornerstone concept in the classical true-score test theory of psychological or educational measurement (e.g.,  Gulliksen, 2013; Lord et al., 1968;).
It is related to measurement error $\varepsilon$, to unexplained or uncontrolled residual variance $\var(\varepsilon)$, 
distinctively inherent to social or behavioral measurements (e.g., Mellenbergh, 1996). 
A normed measure ranging from zero to one, the reliability $\Rel(Y)=\var(\tau)/(\var(\tau)+\var(\varepsilon))$ of an observed test variable $Y=\tau+\varepsilon$
is the proportion of explained or true-score variance $\var(\tau)$, relative to the observed total variance $\var(\tau)+\var(\varepsilon)$. 
Reliability and methods for quantifying reliability, such as Cronbach's alpha (discussed below), 
have been employed in numerous substantial studies (e.g.,  Cortina, 1993; Peterson, 1994; Hogan et al., 2000). 
Reliability is an essential quality criterion required for a ``good'' psychological or educational test, whereby it represents the extent to which a test in repeated independent measurements under same conditions yields comparable test results. 
However, independent test repetitions are not possible, or a test may only be administered once. Thus, reliability is an unobserved or unknown parameter that has to be estimated from empirical data.

Research has examined various coefficients for reliability estimation, see, e.g., the critical discussion published in Psychometrika centered around the paper by  Sijtsma (2009a), with
reactions by Bentler (2009), Green \& Yang (2009a, 2009b), Revelle \& Zinbarg (2009), and Sijtsma (2009b) or also Ten Berge \& So\u{c}an (2004). Thereof, R implementations of the latent class 
reliability coefficient  (van der Ark et al., 2011) or the MS statistic  (Sijtsma \& Molenaar, 1987; Molenaar \& Sijtsma, 1988)  as reliability measures 
can be found in van der Ark (2012). Despite the fact that methodologically superior reliability coefficient exists, one popular and \textit{the} most
widely used method for estimating the reliability of a test or composite score is the coefficient alpha by Cronbach (1951), generally, a lower bound to the reliability -- see also Nunnally \& Bernstein (1978) or 
Furr \& Bacharach (2013).

Compared to such sophisticated methods as the latent class reliability coefficient or the MS statistic, Cronbach's alpha is simple to compute, since the only requirement 
to calculate coefficient alpha is the corresponding covariance matrix. We, thus, exemplify the proposed techniques of this paper based on coefficient alpha, but also 
dicuss applications of the proposed resampling machinery for other reliability measures.

To introduce the coefficient and the corresponding model we consider a test or measurement instrument consisting of $k$ items. 
The observed responses of $N$ examinees, i.e. of $N$ independent and identically distributed repetitions of this test, are denoted by $X_{r1},\dots,X_{rk},$ $1\leq r\leq N,$ and for each examinee the observations are combined into a response vector 
$\bX_r = (X_{r1},\dots,X_{rk})'$. The test or composite score variable of the measurement instrument is the sum $S=\sum_{i=1}^k \overline{X}_{\cdot i}$, where $\overline{X}_{\cdot i}=\frac{1}{N}\sum_{r=1}^N X_{ri}$. It is assumed that the response vectors $\bX_r$ are centered, i.e. $\Ex(\bX_1)=\bnull$, 
and possess a non-zero covariance matrix $\var(\bX_1) = \bSigma$. Then coefficient alpha is defined as  (Cronbach, 1951; Guttman, 1945)
\begin{equation}\label{Cronbachalpha}
   \alpha_C 
  = \alpha_C(\bSigma) 
  = \frac{k}{k-1} \left(1 - \frac{\tr(\bSigma)}{\kone'\bSigma\kone }\right)
  = \frac{k}{k-1} \left(1 - \frac{\sum\limits_{i=1}^k \var(X_{1i})}{\sum\limits_{i,j=1}^k \Cov(X_{1i},X_{1j})}\right),
\end{equation}
where $\kone=(1,\dots,1)'$ denotes the $k$-dimensional vector of ones. 
In the basic classical test theory additive error model  (e.g., Lord et al., 1968), this internal consistency coefficient $\alpha_C$ is only a lower bound for the reliability of the composite score variable $S$. That is, if the error variables
of the test variables are uncorrelated, $\alpha_C \leq \Rel(S)$. However, under the more restrictive model of essentially $\tau$-equivalent variables with uncorrelated residuals  (e.g., Novick \& Lewis, 1966), $\alpha_C$ is equal to 
the reliability of the test score variables $S$, i.e. $\alpha_C=\Rel(S)$. These are restrictive assumptions from a practical point of view. Still, in applications, the coefficient 
alpha is widely used even if these assumptions are not met. Therefore, we have chosen $\alpha_C$ to demonstrate the usefulness of the subsequently presented resampling approaches but also outline extensions to other reliability measures 
in Section \ref{sec:furtherapp} below.
Typically, $\alpha_C$ is estimated by 
\begin{equation}\label{estalpha}
 \widehat{\alpha}_C= \alpha_C(\widehat{\bSigma}) = \alpha_C\left(\frac{1}{N-1} \sum_{i=1}^k (\bX_i - \bmX)(\bX_i - \bmX)'\right),
\end{equation}
where $\bSigma$ in Equation (\ref{Cronbachalpha}) is replaced by the empirical covariance matrix $\widehat{\bSigma}$. 
 Existing inference methods for constructing confidence intervals or statistical tests for $\alpha_C$ in the above described one-sample setting are mostly based on asymptotic results.
In particular,  van Zyl et al. (2000) were the first to study the asymptotic distribution of $A_N := \sqrt{N}\left(\hat{\alpha}_C-\alpha_C\right)$. Under the assumptions of normality 
$\bs X_r \stackrel{i.i.d.}{\sim} \mathcal{N}(\bs 0, \bs\Sigma)$, they could prove that $A_N$ is asymptotically (as $N \rightarrow \infty$) normally distributed with mean zero and a specific variance $\sigma^2$. Together with a consistent estimate of $\sigma^2$ they 
were able to construct confidence intervals and statistical tests for $\alpha_C$ which are valid for $N \rightarrow \infty$, see also Bonett \& Wright (2015) for a recent application. Extensions of the  van Zyl et al. (2000) result to asymptotically 
distribution free (ADF) non-normal models are stated in the fundamental papers by 
Yuan et al. (2003) and  Maydeu-Olivares et al. (2007). Only assuming finite eight moments 
they also obtained that $A_N$ is asymptotically ($N \rightarrow \infty$) normal with mean zero and a 
slightly more complicated limit variance $\tilde{\sigma}^2$. Again asymptotically exact confidence intervals for $\alpha_C$ were determined based on a consistent estimate of $\tilde{\sigma}^2$ (Maydeu-Olivares et al., 2007), nonparametric 
bootstrap techniques (Yuan et al., 2003) or even parametric bootstrap methods  (Padilla et al., 2012). From the simulation study in Padilla et al. (2012) the latter parametric bootstrap procedure 
seemed to be the method of choice. Also, a permutation-type approach for the one-sample case was discussed (see Prelog et al., 2009). 

In comparison to the extensive treatment of the one-sample case, the two-sample case for comparing the reliability of two different subgroups or even two different questionnaires has gained less attention. 
Exceptions are given by the approximate testing procedures of Maydeu-Olivares et al. (2007) as well as Bonett \& Wright (2015) and Bonett (2003), where the latter also contains a certain bootstrap proposal,
where bootstrap samples are drawn separately for each group. This 'within' bootstrap is by far not the best choice in comparison to other resampling methods, 
see e.g. the simulation study in  Konietschke \& Pauly (2014) for a different null hypothesis of interest. 
Alternative resampling procedures do exist that ensure a much better control of the error of first kind (type I error). Permutation methods are particularly suited.
Remark, that such statistical tests are well known for keeping the prescribed level finitely exact if the data is exchangeable (i.e. the joint distribution of the pooled sample does not change under arbitrary permutations of the group status) 
under the null.
However, it is a misapprehension that permutation tests are always valid inference procedure for the larger null hypothesis of interest formulated in terms of parameters or effect measures, 
see e.g. Bradbury's (1987) discussion on the paper by Still \& White (1981). 
Nevertheless, it is less well known that a proper choice of the test statistic may lead to permutation tests that are still 
(at least asymptotically) exact if the data are not exchangeable. 
Examples are given by permutation tests for comparing means  (Janssen, 1997), variances  (Pauly, 2011),
correlations (Omelka \& Pauly, 2012)  or even more complex functionals and designs  (Chung \& Romano, 2013). These authors showed that permutation tests remain valid for larger null hypotheses. 
Recently, Pauly et al. (2015) as well as Umlauft et al. (2017) have applied similar ideas for constructing asymptotically exact permutation tests for general factorial designs. 
In light of these findings, we will adopt the modified permutation approach in this paper to construct a permutation-based testing procedure for comparing 
the population alpha coefficients of two (or even multiple) independent samples of equal dimensions. This means that we have the same number of items in both groups. In case that an unequal number of items/dimensions between the groups is present, 
we will additionally propose an asymptotic model-based bootstrap extending the results of
Padilla et al. (2012). In any case, we provide the theoretical aspects of all resampling methods considered for inference regarding coefficient alpha.
As a byproduct, this also gives a theoretical justification of the 
procedure introduced in Padilla et al. (2012).

The paper is organized as follows. In Section \ref{sec:mod}, the statistical model for the two-sample theory and the resulting asymptotic statistical test are described. A short introduction to two resampling methods is considered in Section \ref{sec:res}. 
First, the permutation approach is described and afterwards, a parametric bootstrap procedure is applied to the two-sample model. Moreover, extensions to other models (one-way layout and paired two-sample designs) are given at the 
end of this section. Section \ref{sec:furtherapp} discusses extensions to other reliability 
measures. In Section \ref{sec:sim}, the usefulness of the resampling-based inference methods are illustrated in a simulation study. 
An application to empirical data of the procedures introduced in this work is given in Section \ref{sec:app}. In Section \ref{sec:dis}, we conclude with a summary, as well as with related remarks on further research.

\section{Statistical Model}\label{sec:mod}

We first explain how the known normal model procedures can be extended to more general models. In particular, we adapt the non-normal asymptotic distribution free (ADF) setting 
described in  Maydeu-Olivares et al. (2007) to our two-sample case by 
considering two groups of independent zero-mean random vectors 
\begin{equation} \label{model}
  \bX_1,\dots, \bX_{n_1} \quad \mbox{ and }\quad  \bX_{n_1+1},\dots, \bX_{N}.
\end{equation}
Here $\bX_r = (X_{r1},\dots,X_{rk_1})'$ for the first sample with $1\leq r\leq n_1$ examinees and
 $\bX_{n_1+s} = (X_{s1},\dots,X_{sk_2})'$ for the second sample with 
 $1\leq s\leq n_2=N-n_1$ 
 and fixed item numbers $k_1,k_2$. In this ADF framework it is only required, that the random vectors are independent and identically distributed in each group with finite eight-order moments ($\Ex(\|\bX_1\|^8)+\Ex(\|\bX_N\|^8)<\infty$)
 and arbitrary covariance matrices $\bSigma_1=\Cov(\bX_1)$ and $\bSigma_2=\Cov(\bX_N)$. Note, that these assumptions are weaker than those given in Kuijpers, Ark, \& Croon (2013) who proposed a different approach based on marignal models.
 
Denote the population alpha coefficients by 
$\alpha_{C,1}= \alpha_C(\bSigma_1)$ and $\alpha_{C,2}= \alpha_C(\bSigma_2)$, respectively, where $\alpha_C$ is as in Equation (\ref{Cronbachalpha}). 
Now, the null hypothesis of interest is given by 
 $H_0 : \alpha_{C,1} = \alpha_{C,2},$ 
which we like to test against one-sided $H_{11} : \alpha_{C,1} > \alpha_{C,2}$ (testing for superiority) or two-sided alternatives $H_{12} : \alpha_{C,1} \neq \alpha_{C,2}$. 
Denoting the empirical covariance matrices of the two-samples by $\widehat{\bSigma}_1$ and $\widehat{\bSigma}_2$, 
respectively,  the population alpha coefficients can be estimated consistently (as $\min(n_1,n_2)\to\infty$) by 
$\widehat{\alpha}_{C,1}= \alpha_C(\widehat{\bSigma}_1)$ and $\widehat{\alpha}_{C,2}= \alpha_C(\widehat{\bSigma}_2)$ since $\alpha_C(\cdot)$ 
is a smooth function. 
Hence, a first naive idea would be to base the statistical test on the following statistic
\begin{equation*}M_n=\sqrt{\frac{n_1n_2}{N}} \left(\widehat{\alpha}_{C,1}-\widehat{\alpha}_{C,2}\right)\end{equation*}
in the one sided case and on $|M_n|$ in the two-sided case. The results from  Maydeu-Olivares et al. (2007) imply that
$\sqrt{n_i}(\widehat{\alpha}_{C,i} - \alpha_{C,i})$ is asymptotically (as $n_i\to\infty$) 
normally distributed with mean zero and variance $\tilde{\sigma}_i^2$ (assumed to be positive) for both the cases $i=1,2$, see the supplementary material to this paper for the explicit form of $\tilde{\sigma}_i^2$. 
In large samples with $\frac{n_1}{N}\rightarrow \kappa \in(0,1)$, the asymptotic null distribution of $M_n$ is thus given by a standard normal distribution with mean zero and variance $\tilde{\sigma}^2=(1-\kappa)\tilde{\sigma}_1^2 + \kappa\tilde{\sigma}_2^2$
since both groups are independent.
 Since the limit variance $\tilde{\sigma}^2$ is unknown it would be possible to directly apply a resampling procedure like a bootstrap or permutation approach to this result.  However, it turns out that for obtaining an adequate permutation procedure it is necessary to studentize this statistic, 
see e.g. the discussions in  Janssen (1997), Chung \& Romano (2013) or Pauly et al. (2015). 
In particular, a consistent estimator $\widehat{\sigma}^2$ for the limit variance $\tilde{\sigma}^2$ 
of $M_n$ can be adopted from  Maydeu-Olivares et al. (2007) and Maydeu-Olivares et al. (2010), see Equation~\eqref{pooled variance nonormal} of the supplement for its definition. With its help, we can define a studentized test statistic by
\begin{equation}\label{teststatnonnormal}
 T_n =  T_n(\mathbb{X}) = \frac{M_n}{\widehat{\sigma}},
\end{equation}
which is asymptotically equivalent to the ADF test statistic considered in Maydeu-Olivares et al. (2007, Equation 4). 
It follows that $ T_n $ is asymptotically standard normal and the corresponding one-sided asymptotic exact level-$\alpha$\footnote{$\alpha$ without any indices denotes the significance level of 
the corresponding statistical test, not Cronbach's alpha.}-test 
$\varphi_n=\mathbbm{1}\{T_n>z_{1-\alpha}\}$ compares the test statistic $T_n$ with
the $(1-\alpha)$-quantile of a standard normal distribution $z_{1-\alpha}$ and rejects $H_0$ for large values of $T_n$. Since the finite sample properties of this inference method is in general not desirable (see Section \ref{sec:sim}), different 
resampling principles are proposed to improve its small sample behaviour.

\section{Proposed Resampling Approaches}\label{sec:res}

\subsection{Permutation Procedure}\label{sec:perm}
Let us shortly recall the general permutation idea assuming an equal number of items. Instead of basing the statistical testing method on the asymptotic results for $T_n$, i.e. 
choosing quantiles from a standard normal distribution as critical values, the critical values are obtained as quantiles from the corresponding conditional
permutation distribution of $T_n$ given $\bs X$. The reason for this is twofold. 
First, this approach yields an inference method that even keeps the prescribed level for finite sample sizes exactly if the joint distribution of 
$\bX_1,\dots,\bX_N$ is invariant under random permutation of the group status (i.e. the vector is exchangeable). 
In the special case of multivariate normality, this means that the covariance matrices in both groups are equal $\bSigma_1=\bSigma_2$. 
Second, if this is not the case, the permutation distribution of $T_n$ is data-dependent 
and hence intuitively mimics the unknown null distribution of $T_n$ at finite sample size better than the asymptotic normal approximation which results in  
more adequate test decisions. A detailed theoretic explanation for this procedure is given in the  supplement to this paper.

Now, this resampling approach will be formalized:
Given the observed responses, let 
$\mathbb{X}^\pi = (\bX_{\pi(1)},\dots,\bX_{\pi(N)})$ denote a random permutation of all $N$ data vectors
$\mathbb{X}=(\bX_1,\dots,\bX_N)$, i.e. the group status in $\mathbb{X}^\pi$ is randomly given (without replacement). 
Here, $k_1=k_2=k$ is assumed, $\pi$ is a random permutation that is independent of the responses and uniformly distributed on the set of all permutations of the numbers $1,\dots,N$.  
Note, that we only permute the vectors (i.e. test repetitions) and not all components (i.e. item responses).
Assuming $k_1=k_2$, the values of $M_n$ and its variance estimator $\hat{\sigma}^2$ are calculated from 
the permuted observations $\mathbb{X}^\pi$ to obtain the permutation version of the test statistic
\begin{equation}\label{equ:tildeTn}
  T_n^\pi =  T_n(\mathbb{X}^\pi) = \frac{M_n(\mathbb{X}^\pi)}{\tilde{\sigma}(\mathbb{X}^\pi)}.
\end{equation}

Let now $c_n^{\pi}(\alpha)$ denote the $(1-\alpha)$-quantile 

of the conditional permutation distribution of $T_n^\pi$ given $\bX_1,\dots,\bX_N$, i.e. of 
\begin{equation}\label{permdist}
  x\mapsto \frac{1}{N!} \sum_{\pi} \mathbbm{1}\{T_n(\mathbb{X}^\pi)\leq x\},
\end{equation}
where the summation is over all $N!$ possible permutations. Note, that due to symmetry in the test statistic only $\binom{N}{n_1}$ different summands have to be calculated in practice.
For larger $N$, however, their calculation is computationally too expensive and the permutation distribution function (\ref{permdist}) is usually approximated via Monte-Carlo methods, see the algorithm below.
Then our proposed studentized permutation test is given by 
\begin{equation}\label{permtest}
\psi_n = 
\mathbbm{1}\{T_n > {c}_n^{\pi}(\alpha)\}  + \gamma_n^\pi(\alpha) \mathbbm{1}\{T_n = c_n^{\pi}(\alpha)\}
\end{equation}
in the one-sided case. 
 Note, that the randomization $\gamma_n^\pi(\alpha)$ can be omitted for large sample sizes since the test statistic is asymptotically continuously distributed. 
 However, since no specific 
 assumption of the shape of the distribution is made (rather than some moment assumptions) it is in general needed to gain finite exactness of the 
 statistical test if the distribution of the data vector is 
 invariant under permutation of the group status. For example in the scale-model 
 \begin{equation}\label{model2}
  \bX_i= \bSigma_j^{1/2}\bep_i, \qquad\text{ where } j=1 \text{ for } 1\leq i\leq n_1 \text{ and } j=2 \text{ otherwise,}
 \end{equation}
 with i.i.d. random vectors $\bep_i$  and covariance matrix 
 given by the identity $\bs I_k$, this invariance property is fulfilled 
 iff $\bSigma_1=\bSigma_2$.

In the supplementary material it is shown that $\psi_n$ is an asymptotically exact testing procedure, i.e. its error of the first kind is approximate $\alpha$ for large sample sizes, provided that 
$\frac{n_1}{N} - \kappa= O(n_1^{-1/2}).$ 
As mentioned above it is even an exact level $\alpha$ testing for finite sample sizes 
if the data is exchangeable.
Finally, we additionally state the algorithm for the computation of the p-value for the overall sample size $N$ in the two-sided case:
\begin{enumerate}
\item Given the data $\mathbb{X}=(\bX_1,\dots,\bX_N)$, compute the studentized test statistic $T_n=T_n(\mathbb{X})$ as given in Equation (\ref{teststatnonnormal}).
\item Obtain $\mathbb{X}^\pi=(\bX_{\pi(1)},\dots,\bX_{\pi(N)})$ by randomly permuting the data vectors.
\item Calculate the permuted version of the test statistic $T_n^\pi=T_n(\mathbb{X}^\pi)$ as in Equation \eqref{equ:tildeTn}.
\item Repeat the steps 2.-3. $B$ (e.g. $10^3$ or $10^4$) times and save the values $T_n^\pi$ in  $A_1,\ldots,A_{B}$. 
\item Estimate the two-sided p-value by $$p = \min\{2p_1,2-2p_1\}, \; \text{where}\; 
p_1=\frac{1}{B}\sum_{\ell=1}^{B} \mathbbm{1}\{T_n\leq A_\ell\}. $$
\end{enumerate}
In the one-sided case the p-value is estimated by $1-p_1$.
Note, that we have omitted the randomization $\gamma_n^\pi(\alpha)$ of the two-sided statistical test for ease of convenience. 
Also a corresponding permutation-based two-sided confidence interval can be calculated:
\[\left[(\hat{\alpha}_{C,1}-\hat{\alpha}_{C,2})\pm \frac{c^\pi_n (\nicefrac{\alpha}{2})}{\sqrt{\frac{n_1n_2}{N}}}\cdot \hat{\sigma}\right].\]
It possesses asymptotic coverage probability of level $1-\alpha$ for the 
unknown Cronbach's alpha differences $\alpha_{C,1}-\alpha_{C,2}$. 
Since  the proposed permutation procedure is only applicable for an equal number of test items $k_1=k_2$ 
we additionally study a parametric bootstrap technique that is even applicable in case of possibly different item sizes.

\subsection{Parametric Bootstrap Procedure}\label{sec:boot}
We now study the parametric bootstrap procedure, which  Padilla et al. (2012) have applied with regard to coefficient alpha in the one-sample situation, 
see also Konietschke et al. (2015) for a recent application in the MANOVA context.
For observed sample covariances $\widehat{\bSigma}_1$ and $\widehat{\bSigma}_2$ the resampling mechanisms is given by generating independent 
bootstrap variables
\begin{equation*}
  \bX_{1}^\star, \dots, \bX_{n_1}^\star \stackrel{i.i.d.}{\sim} \mathcal{N}_{k_1}(\bnull, \widehat{\bSigma}_1) \quad \mbox{ and } \quad
  \bX_{n_1+1}^\star, \dots, \bX_{N}^\star \stackrel{i.i.d.}{\sim} \mathcal{N}_{k_2}(\bnull, \widehat{\bSigma}_2).
\end{equation*}

With this we calculate a parametric bootstrap version $T_n^\star=T_n(\bX_{1}^\star,\dots,\bX_{N}^\star)$ 
of our studentized test statistic which is used to approximate the unknown distribution of $T_n$ in Equation \eqref{teststatnonnormal}. 
In this case, it follows from a pointwise application of the multivariate CLT that, given the observed responses, 
the distribution of $T_n^\star$ is asymptotically standard normal in probability, see the supplement for the derivation. 
Hence, the one-sided parametric bootstrap test 
\begin{equation}\label{pbtest}
\psi_n^\star = 
\mathbbm{1}\{T_n > c_n^{\star}(\alpha)\} 
\end{equation}
is also asymptotically exact, where $c_n^{\star}(\alpha)$ is the $(1-\alpha)$-quantile of the conditional parametric bootstrap distribution function of $T_n^\star$ given the data $\bX_1, \ldots, \bX_N$. The notion parametric bootstrap may be 
misleading, (since it is also valid in the general ADF framework as shown in the supplement) and asymptotic model-based bootstrap may be a more appropriate term. 
Since $\mathbb{P}\left(|T_n^\star|\leq c^\star_n (\nicefrac{\alpha}{2})\right)\approx 1-\alpha$ for large sample sizes, 
we additionally obtain a corresponding two-sided confidence interval of approximate level $1-\alpha$ given by
\begin{equation*}
 \left[(\hat{\alpha}_{C,1}-\hat{\alpha}_{C,2})\pm \frac{c^\star_n (\nicefrac{\alpha}{2})}{\sqrt{\frac{n_1n_2}{N}}}\cdot \hat{\sigma}\right].
\end{equation*}
Due to the similarity of the permutation and bootstrap test, the confidence intervals are nearly the same. Only the quantiles differ.
 To compare the finite sample performance of $\psi_n^\star$  with that of the permutation test $\psi_n$ and the asymptotic benchmark $\varphi_n$ 
a simulation study is conducted in Section~\ref{sec:sim}.
We like to stress that \textsc{R} code for carrying out both resampling (permutation and parametric bootstrap) procedures is given in the supplementary material to this paper.

\subsection{Remarks on the ADF Assumption and Extensions to One-Way Layouts}\label{sec:ext}

\noindent{\bf On the Distributional Assumptions.} 
We have consciously chosen to work under a general ADF framework since multivariate normality is a rather strong assumption that is 
usually violated for practical data at hand, see e.g. the discussion in  Konietschke et al. (2015). 
This is especially the case when confronted with ties in the data and / or small sample sizes. 
For completeness, however, we like to point out that similar but computationally simpler resampling procedures can be derived by just changing the 
consistent variance estimator in the definition of $T_n$ to the more simple estimate of van Zyl et al. (2000) in the normal case, leading to a more simple test statistic, say $\tilde{T}_n$. The resulting asymptotic test would be 
related to the two-sample test considered in Bonett \& Wright (2015). 
It is then straightforward to prove that parametric bootstrap and permutation procedures based on 
$\tilde{T}_n$ are also valid under the normality assumption. 
In fact, our proposed parametric bootstrap procedure originally stems from this parametric model. The motivation to apply it 
with a slightly different covariance estimator also in the ADF case is due to the multivariate central limit theorem, see also the explanation in Konietschke et al. (2015). \\

\noindent{\bf Extensions to Multiple Samples.} 
The above inference procedures for two independent groups can also be extended to compare the Cronbach coefficient $\alpha_C$ from $K$ 
independent samples. 
In particular, denoting the corresponding population alpha coefficients in group $j$ by $\alpha_{C,j}$, ($j=1,\dots,K$) 
this leads to the null hypothesis $H_0^K: \alpha_1 = \ldots = \alpha_K$ that has also been considered in  Kim \& Feldt (2008). 
To fix notation, let $\hat{\sigma}_j^2$ be the consistent variance estimator of Maydeu-Olivares et al. (2007) in group $j$ and denote 
 the sample size of the $j$-th group by $n_j$ and the total sample size by $N=\sum_{j=1}^K n_j$.
Writing $\hat{\pmb \alpha}_{C} = (\hat{\alpha}_{C,1},\dots, \hat{\alpha}_{C,K})'$ and $\hat{\pmb \Sigma} = \diag(\frac{N}{n_j} \hat{\sigma}_j^2, \; j=1,\dots,K)$, a suitable test statistic for 
$H_0^K$ is given by 
$$
Q_N = N \hat{\pmb \alpha}_{C}' {\pmb H}_K ({\pmb H}_K \hat{\pmb \Sigma} {\pmb H}_K)^+ {\pmb H}_K \hat{\pmb \alpha}_{C}.
$$
Here, ${\pmb H}_K = {\pmb I}_K - K^{-1} {\pmb 1}_K{\pmb 1}_K'$ and $(\cdot)^+$ denotes the Moore-Penrose inverse. 
Combining the techniques from the supplementary material with results from  Chung \& Romano (2013, Theorem 3.1) and Konietschke et al. (2015) we can derive valid permutation (for equal numbers of group items) and bootstrap procedures for $H_0^K$ by comparing 
$Q_N$ with critical values taken from its corresponding resampling version. For $K=2$ the $Q_N$-based inference procedures simplify to the bootstrap and permutation tests in $T_n$ from above. 
In practice, this may subsequently lead to a hierarchical multiple testing problem: 
After rejecting $H_0^K$ one may test all $\binom{K}{2}$ pairwise two-sample hypotheses 
$H_0^{(i,j)}: \alpha_i = \alpha_j$, $1\leq i < j\leq K,$ by means of the tests from Section~\ref{sec:perm}-\ref{sec:boot}; 
possibly adjusted for multiplicity.\\

\noindent{\bf Extensions to paired designs.}

Moreover, the parametric bootstrap procedure is also applicable for paired two-sample designs, where data is given by the
independent and identically distributed random vectors
\[\bs{X}_i=(\bs{X}'_{1,i}, \bs{X}'_{2,i})', i=1,\ldots,N.\] Here $\bs{X}_{1,i}$ and $\bs{X}_{2,i}$ contain all observations of individual $i$ for treatment / time point 1 and 2, respectively. Let $\bSigma_1=\Cov(\bX_{1,1})$ and $\bSigma_2=\Cov(\bX_{2,1})$ denote the corresponding covariance matrices, whereas 
$\bSigma_{12}=\Cov(\bX_{1,1}, \bX_{2,1})$ describes the covariance structure of the pairs. Altogether, 
this leads to a covariance matrix $\bSigma=\Cov(\bX_1)=\begin{pmatrix}\bSigma_1 & \bSigma_{12} \\
                                                                                                                                        \bSigma_{12} & \bSigma_2
                                                                                                                                       \end{pmatrix}$.
The null hypothesis for paired data is again given as $H_0^{\text{pair}}: \alpha_{C,1} =\alpha_{C,2}$, where $\alpha_{C,1}=\alpha_C(\bSigma_1)$ is the corresponding alpha coefficient for the first and $\alpha_{C,2}=\alpha_C(\bSigma_2)$ denotes Cronbach's alpha for the second time point, respectively. 
Since the sample covariance fulfills a central limit theorem (i.e. $\sqrt{N}(\widehat{\bSigma}-\bSigma)$ is asymptotically normal) it thus follows from an application of the delta-method and Slutzky's Lemma that $P_N = \sqrt{\frac{n_1 n_2}{N}} \left(\alpha_C(\widehat{\bSigma}_1)-\alpha_C(\widehat{\bSigma}_2)\right)$ is asymptotically normally distributed  
with some specific variance $b^2$. Assuming $b^2>0$ and denoting the obvious plug-in estimate as $\hat{b}^2$, we obtain an ADF procedure based on the studentized test statistic ${P_N}/{\hat{b}}=: T_N$ and $z$-quantiles as critical values.              
In addition, we can apply a modification of the parametric bootstrap procedure from above. Here, the resampled data is given by generating independent bootstrap samples
\[\bs{X}^{\star\star}_1, \ldots, \bs{X}^{\star\star}_N\stackrel{i.i.d.}{\sim}\mathcal{N}(\bs{0}, \widehat{\bSigma}),\]
where $\widehat{\bSigma}=\begin{pmatrix}
                           \widehat{\bSigma}_1 & \widehat{\bSigma}_{12} \\
                           \widehat{\bSigma}_{12} & \widehat{\bSigma}_2
                         \end{pmatrix} $
denotes the observed sample covariance of the whole data matrix $\bs{X}=(\bX_1, \ldots, \bX_N)$. To test the hypothesis given above you again have to calculate the parametric bootstrap version $T^{\star\star}_N=T_N(\bX_1^{\star\star}, \ldots, \bX_N^{\star\star})$ 
of the novel studentized test statistic to compute critical values.  Altogether, this yields an adequate test procedure for paired data.

\section{Other reliability functionals}\label{sec:furtherapp}

In this Section, we sketch that the use of the proposed resampling and inference principles are not limited to Cronbach's alpha coefficient as a reliability
measure. In particular, the applicability of the discussed resampling methods can be extended to certain smooth functionals of covariances. To this end, we discuss some of the reliability measures presented in Revelle \& Zinbarg (2009).
The first category of measures introduced therein are based on the work of Guttman (1945) and lead to the following six $\lambda$-reliability-coefficients:
\begin{center}
 \begin{tabular}{ll}
  $\lambda_1 =1-\frac{\tr(\bSigma)}{\kone'\bSigma\kone}$, & $\lambda_2 =\frac{\kone'\bSigma\kone - \tr(\bSigma) + (\frac{k}{k-1}C_2)^{\frac{1}{2}}}{\kone'\bSigma\kone}$,\\

  $\lambda_3 = 2\left(1-\frac{\kone'\bSigma_A\kone+\kone'\bSigma_B\kone}{\kone'\bSigma\kone}\right)$,  &
  $\lambda_4 = \lambda_1+\frac{2(\bar{C}_2)^{\frac{1}{2}}}{\kone'\bSigma\kone}$,\\
  $\lambda_{5} = \lambda_1+ \frac{k}{k-1}\frac{2(\bar{C}_2)^{\frac{1}{2}}}{\kone'\bSigma\kone}$, &
  $\lambda_6 = 1-\frac{\sum_{t=1}^k e^2_t}{\kone'\bSigma\kone}$,
\end{tabular}
\end{center}
where $C_2=\kone(\bSigma-\diag(\bSigma))^2\kone'$, $\bSigma_A$ and $\bSigma_B$ are obtained by splitting $\bSigma$ into two parts (no matter how the test is splitted) and $e^2_t$ are the variance of the errors. Moreover, $\bar{C}_2$ denotes the maximal value of $C_{2t}$, the sum of squares of 
the covariances of item $t$. Regarding these six measures, it is assumed that the covariances between the items 
represent the true covariance, whereas the variance matrix $\bSigma$ reflects an unknown sum of true ($\bSigma_t$) and some error variances ($\bSigma_e$), i.e. $\bSigma=\bSigma_t+\bSigma_e$. 

The second category of reliability measures given in Revelle \& Zinbarg (2009) are based on a decomposition of the variance into four parts: a general factor $\bs{g}$, a  
group factor $\bs{f}$, a specific factor $\bs{s}$ which is unique to each item, and a random error $\bs{e}$ (see McDonald, 1978, 1999). All these factors
are combined to obtain the following model for the vector $\bs x$ of observed scores  in the $k$ scale items
\[\bs x=\bs{cg}+\bs{Af}+ \bs{Ds}+\bs{e},\]
where $\bs c$ is a vector of general factor loadings, $\bs A$ a matrix of group factor loadings 
and $\bs D$ a diagonal matrix of factor loadings on the item specific factors  (see Zinbarg et al., 2005).
McDonald (1978, 1999) propose the following reliability measures:
\begin{equation*}
  \lambda_7 = \frac{\kone'\bs{cc}'\kone+\kone'\bs{AA}'\kone}{\kone'\bSigma\kone} \;\; \text{ and } \;\;
  \lambda_8 = \frac{\kone'\bs{cc}'\kone}{\kone'\bSigma\kone}.
\end{equation*}

All of these measures have in common that they are smooth functionals of the underlying covariance matrix $\bSigma$ (see the supplement for details). Thus, it follows from the multivariate delta method and the CLT of the sample covariances that 
$M_{n,i}=\sqrt{\frac{n_1 n_2}{N}}\left(\lambda_i(\widehat{\bSigma}_1)-\lambda_i(\widehat{\bSigma}_2)\right),$ $i=1,\dots,8,$ is asymptotically normal under $H_{0,i}: \lambda_i(\bSigma_1)=\lambda_i(\bSigma_2)$ with mean zero and some specific variance $\sigma^2_{\lambda_i}$
under the model assumptions stated in Section \ref{sec:mod}. Again consistent plug-in estimates $\hat{\sigma}^2_{\lambda_i}$ for $\sigma^2_{\lambda_i}$ are obtained easily leading to ADF procedures for testing $H_{0,i}$ based on the test statistic $M_{n,i}$,
the variance estimates $\hat{\sigma}^2_{\lambda_i}$ and a $z$-quantile as a critical value. Moreover, carefully
checking the proof of the parametric bootstrap approach the corresponding results directly carry over. In particular, due to a conditional multivariate central limit theorem for the parametric bootstrap versions $\widehat{\bSigma}_1^\star$ and $\widehat{\bSigma}_2^\star$
of the sample covariances, the validity of its applicability can be proven by means of the delta method for the bootstrap.
For the permutation procedure, it is additionally needed that the corresponding empirical estimators (based on the sample covariance) are asymptotically linear. Since 
all the reliability measures described above are polynomials of the covariances, these functionals are smooth and therefore differentiable in $\bSigma$. Applying Taylor eventually leads to the desired results.

\section{Simulation Study}\label{sec:sim}

This section investigates the properties of the proposed techniques within a simulation study. To make our results comparable to the simulation results of Maydeu-Olivares et al. (2007), 
ordinal data were generated following an algorithm by Muthén \& Kaplan (1985, 1992), where ordinal variables are assumed to be resulting from discretized continuous variables based on thresholds. Note, that
the differences of the item scores are assumed to be meaningful. Simulations regarding continuous data are provided in the supplementary material.
The simulations are conducted with the help of \textsc{R} computing environment, version 3.2.0  (R Core Team, 2016).

\subsection{Design} 
Multivariate normal data with mean zero and correlation matrix $\bs{P}$ were generated. To discretize these data, 
a vector of thresholds $\bs{\tau}$ was used.
We compared the statistical tests $\psi_n$, $\varphi_n$ and $\psi_n^\star$. 
The first one describes the permutation test, the second one the asymptotic test and the third test uses the 
parametric bootstrap procedure for non-normal models.
Overall, 256 conditions were examined:
\begin{enumerate}
 \item Eight different sample sizes ($n_1$, $n_2$): (10,10), (10,20), (25,25), (25,50), (50,50), (50,75), (75,75), (75,100).
 \item Two different test lengths ($k_1=k_2$): 5 and 20 items.
 \item Eight different correlation matrices: $\bs{P_1}, \ldots, \bs{P_8}$.
 \item Two different choices of thresholds: $\bs{\tau}_1$ and  $\bs{\tau}_2$.
\end{enumerate}
A constant number of five item categories was chosen. The first category was set to zero, thus, the items were scored $0, 1, \ldots, 4$. To investigate the tests' finite sample properties, 
rather small sample sizes were specified. For example, this is useful in school comparison educational studies, where student samples of approximately school class sizes are surveyed and compared. 
The different numbers of items are the shortest and longest lengths typically employed for achievement tests. In case of five items, the eight different correlation matrices were chosen as:
\begin{itemize}
 \item $\bs{P_1} = 0.16 \cdot \bs{J_5} + (1-0.16) \cdot \bs{I_5}$, \quad $\bs{P_2} = 0.36 \cdot \bs{J_5} + (1-0.36) \cdot \bs{I_5}$,
 \item $\bs{P_3} = 0.64 \cdot \bs{J_5} + (1-0.64) \cdot \bs{I_5}$, \quad $\bs{P_4} = \bs{\lambda} \cdot \bs{\lambda}' + \bs{I_5} - \diag(\bs{\lambda} \cdot \bs{\lambda}')$, 
 \item $\bs{P_5} = 0.16 \cdot \bs{J_5} + \diag(0.84, 0.74, 0.64, 0.54, 0.44)$,
 \item $\bs{P_6} = 0.36 \cdot \bs{J_5} + \diag(0.64, 0.54, 0.44, 0.34, 0.24)$,
 \item $\bs{P_7} = 0.64 \cdot \bs{J_5} + \diag(0.36, 0.31, 0.26, 0.21, 0.16)$ and
 \item $\bs{P_8} = \bs{\lambda} \cdot \bs{\lambda}'  - \diag(\bs{\lambda} \cdot \bs{\lambda}') + \diag(1,0.9, 0.8, 0.7, 0.6)$, with $\bs{\lambda}$ = $(0.3, 0.4, 0.5, 0.6, 0.7)'$,
\end{itemize}
where $\bs J_5=\bs 1_5 \bs 1_5'$ denotes the 5-dimensional matrix of ones and $\bs I_5$ the 5-dimensional unit matrix.

\begin{figure}[!ht]
\centering
 \subfigure[vector of thresholds $\tau_1$]{\includegraphics[width=0.3\textwidth]{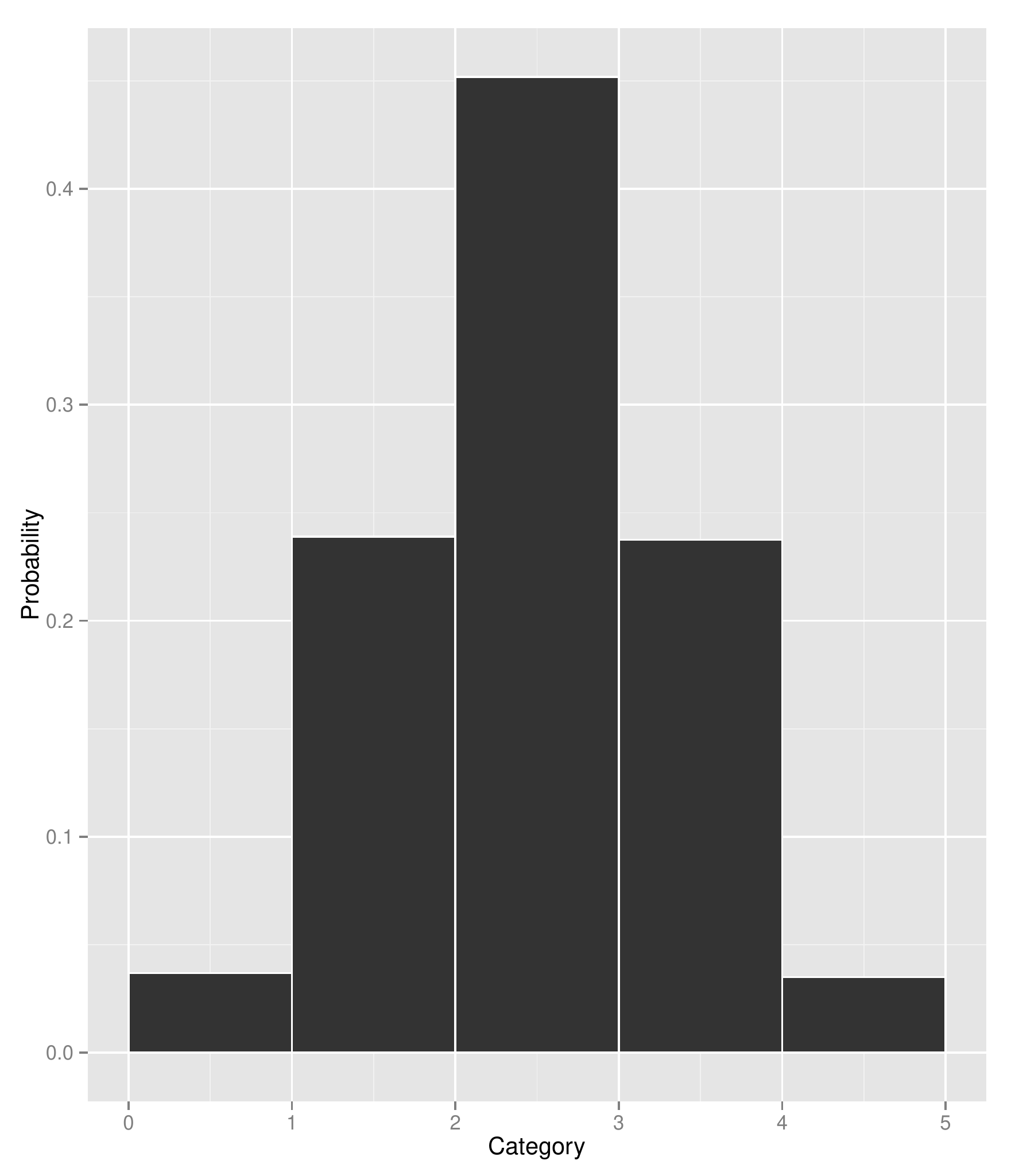}}\qquad
 \subfigure[vector of thresholds $\tau_2$]{\includegraphics[width=0.3\textwidth]{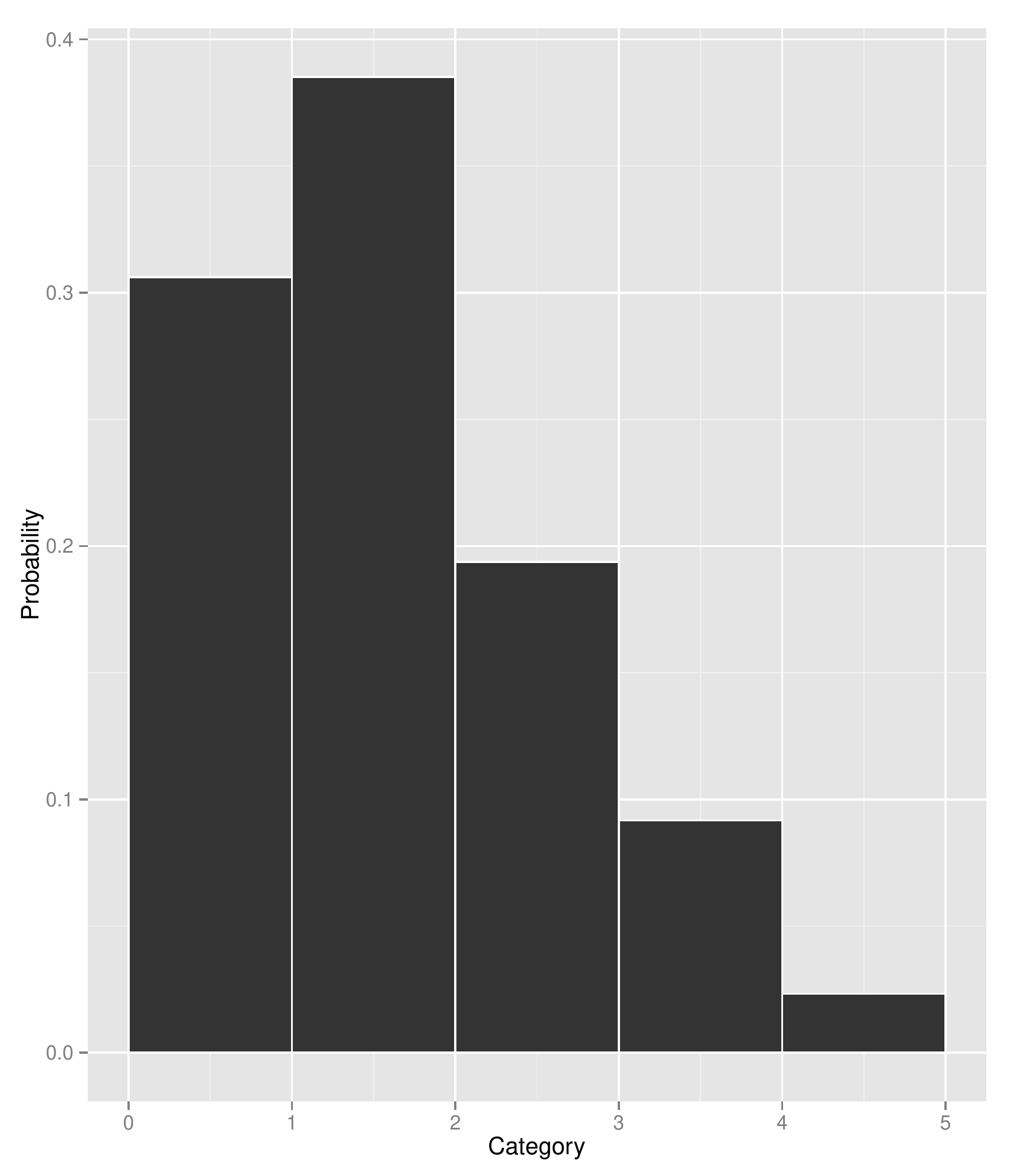}}\\
 \caption{\label{fig:threshold}Histograms of the two different types of items used in this simulation study for threshold $\bs\tau_1$ (a) and threshold $\bs\tau_2$ (b).}
\end{figure}

In case of twenty items the correlation matrices $\bs{P_1}, \bs{P_2}, \bs{P_3}$ are the same compared to the matrices concerning five items, except for the dimensions of $\bs J_{20}$ and $\bs I_{20}$. $\bs{P_4}$ and $\bs{P_8}$ were generated similar to
the correlation matrices for the 5-item case, only
regarding the vector $\bs\lambda=(0.32, 0.34, \ldots, 0.70)' \in \mathbb{R}^{20}$. The last diagonal matrix in the equation of 
$\bs{P_8}$ is generated by using the vector $(0.98, 0.96, \ldots, 0.6)'=(0.98-i\cdot 0.02)_{i=0}^{19}$. The diagonal matrices of $\bs{P_5}, \bs{P_6}, \bs{P_7}$ were chosen as follows: For $\bs{P_5}$ the vector for generating the diagonal matrix is denoted 
by $(0.82-i\cdot 0.02)_{i=0}^{19}$, for $\bs{P_6}$ the vector is given by $(0.62-i\cdot 0.02)_{i=0}^{19}$ and for $\bs{P_7}$ by $(0.35-i\cdot 0.01)_{i=0}^{19}$. The choice of the correlation matrices is based on the simulation studies conducted in 
 Maydeu-Olivares et al. (2007).
One requirement for estimating the true reliability of the test score is the true-score equivalence of items in the population (Lord et al., 1968, Chapter 6). This true-score equivalent model is a model in which the factor loadings 
are equal for all items. This implies that the covariances of the population are all the same, whereas the variances 
are not necessarily equal for all items. Some  of the correlation matrices described beyond ($\bs{P_1}, \bs{P_2}, \bs{P_3}, \bs{P_5}, \bs{P_6}, \bs{P_7}$) follow this assumption and others ($\bs{P_4}$ and $\bs{P_8}$) do not.
The different choices of thresholds ($\bs{\tau}_1, \bs{\tau}_2$) adjusted for the 
skewness and/or kurtosis of the data. In the left panel of Figure \ref{fig:threshold}, the histogram of the data has the form of a normal distribution. In this case, the vector of thresholds $\bs \tau_1$ has entries $(-1.8, -0.6, 0.6, 1.8)$. In the right panel, 
the histogram is shifted to the left. The threshold vector is given by $\bs \tau_2=(-0.4, 0.5, 1.2, 2)$. 

For each of the 256 combinations of sample size, test length, correlation matrix and threshold vector, 
$10,000$ simulation trials were performed, where, for each trial, the results of the statistical tests $\psi_n$ and $\psi_n^\star$ were calculated from a total of $1,000$ permutation or $1,000$ bootstrap samples, respectively. 

\subsection{Results}

\begin{figure}[!ht]
 \centering
\includegraphics[width=0.95\textwidth]{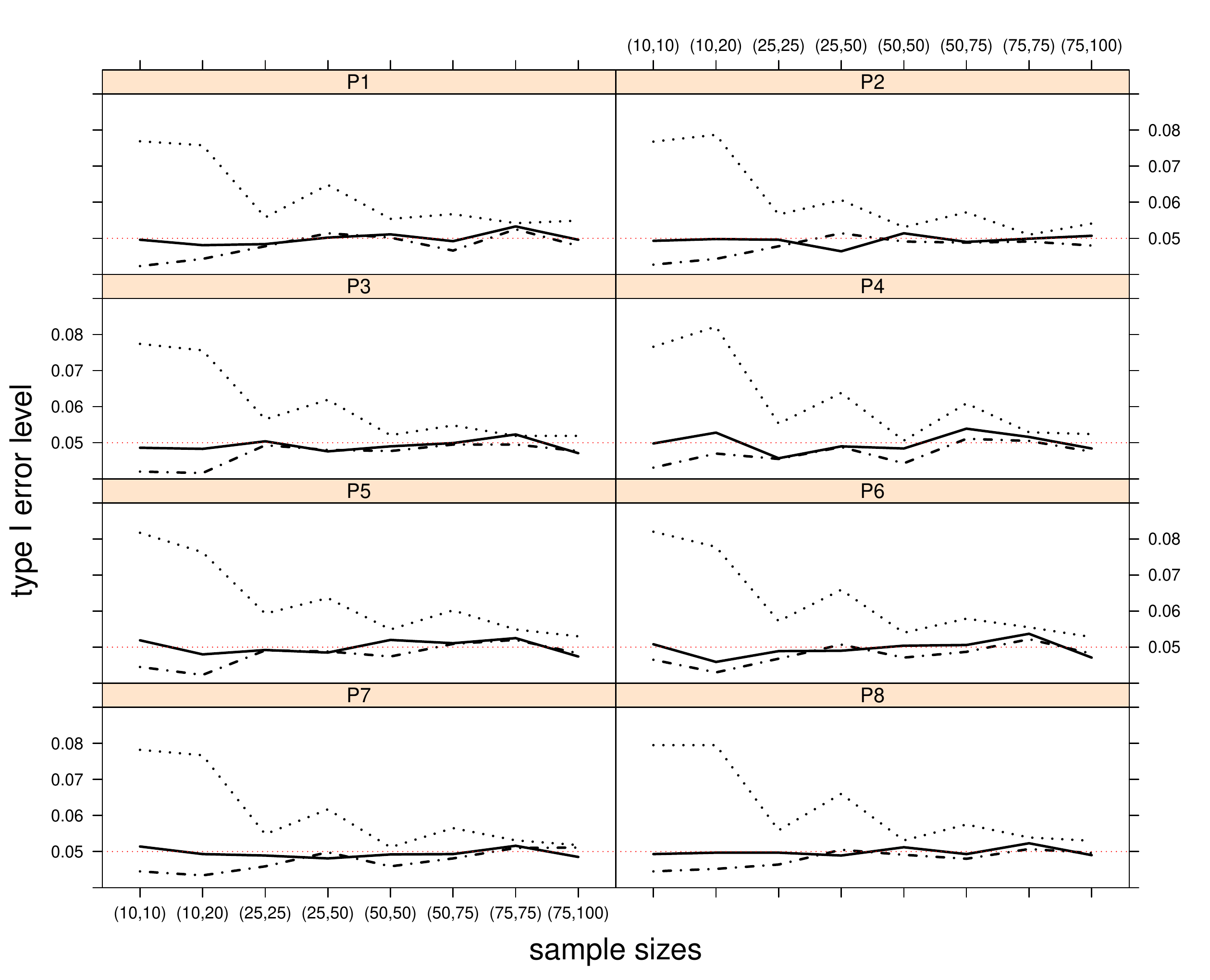}
\caption[kein Strich]{\label{fig:erg5_1}Type I error level ($\alpha=5\%$)  simulation results (y-axis) for threshold $\bs\tau_1$ and 5 items of the permutation test $\psi_n$
  (\rule[0.5ex]{0.7cm}{1pt}), the asymptotic test $\varphi_n$ (\hdashrule[0.5ex]{0.7cm}{1pt}{1pt})  and the bootstrap test $\psi_n^\star$ (\hdashrule[0.5ex]{0.7cm}{1pt}{1pt 1pt 5pt 1pt}) for different sample sizes (x-axis).}
\end{figure}

\begin{figure}[!ht]
 \centering
\includegraphics[width=0.95\textwidth]{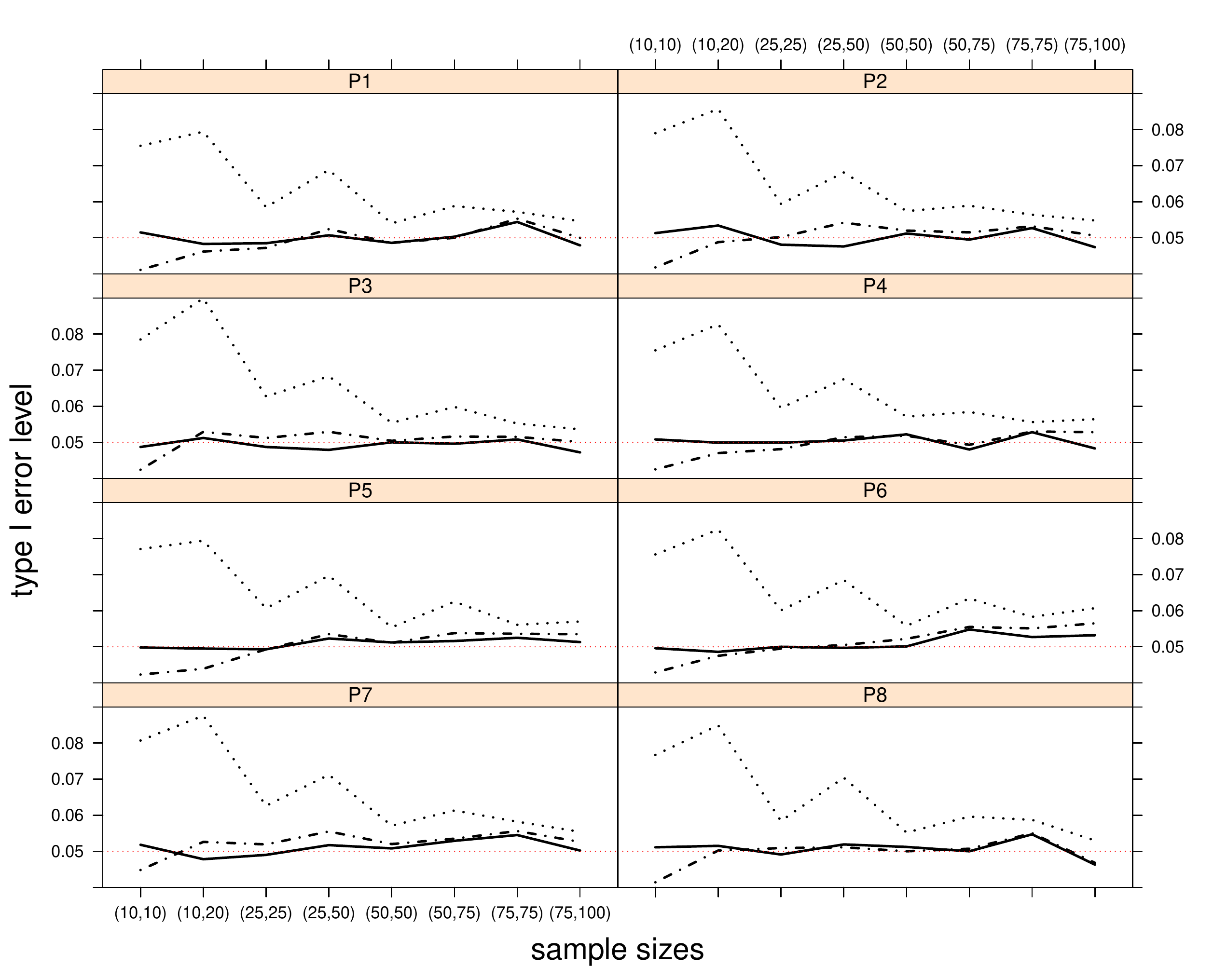}
\caption[kein Strich]{\label{fig:erg5_2}Type I error level ($\alpha=5\%$)  simulation results (y-axis) for threshold $\bs\tau_2$ and 5 items of the permutation test $\psi_n$
  (\rule[0.5ex]{0.7cm}{1pt}), the asymptotic test $\varphi_n$ (\hdashrule[0.5ex]{0.7cm}{1pt}{1pt}) and the bootstrap test $\psi_n^\star$ (\hdashrule[0.5ex]{0.7cm}{1pt}{1pt 1pt 5pt 1pt})
  for different sample sizes (x-axis).}
\end{figure}

\begin{figure}[!ht]
 \centering
\includegraphics[width=0.95\textwidth]{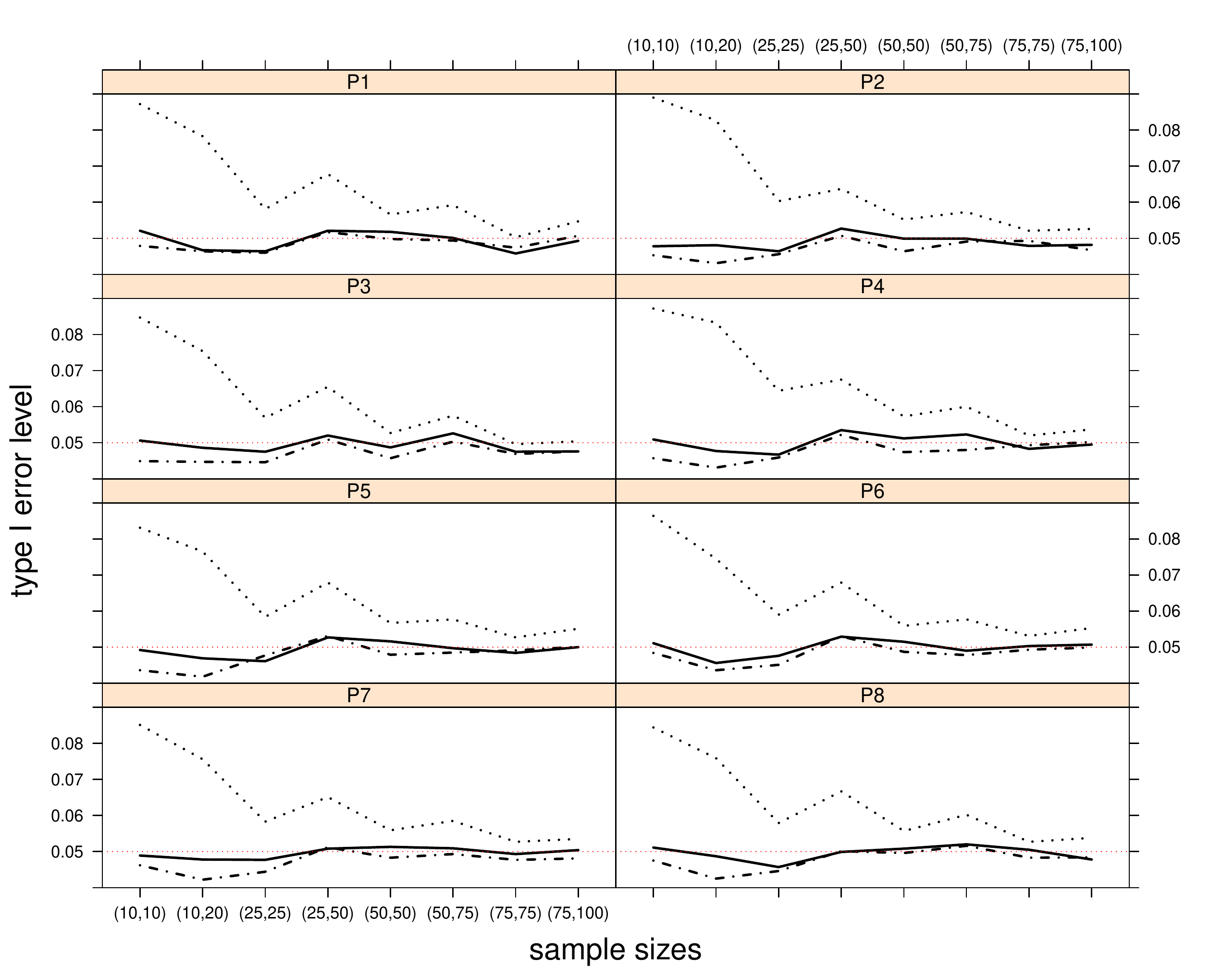}
\caption[kein Strich]{\label{fig:erg20_1}Type I error level ($\alpha=5\%$)  simulation results (y-axis) for threshold $\bs\tau_1$ and 20 items of the permutation test $\psi_n$
  (\rule[0.5ex]{0.7cm}{1pt}), the asymptotic test $\varphi_n$ (\hdashrule[0.5ex]{0.7cm}{1pt}{1pt}) and the bootstrap test $\psi_n^\star$ (\hdashrule[0.5ex]{0.7cm}{1pt}{1pt 1pt 5pt 1pt}) 
  for different sample sizes (x-axis).}
\end{figure}

\begin{figure}[!ht]
 \centering
\includegraphics[width=0.95\textwidth]{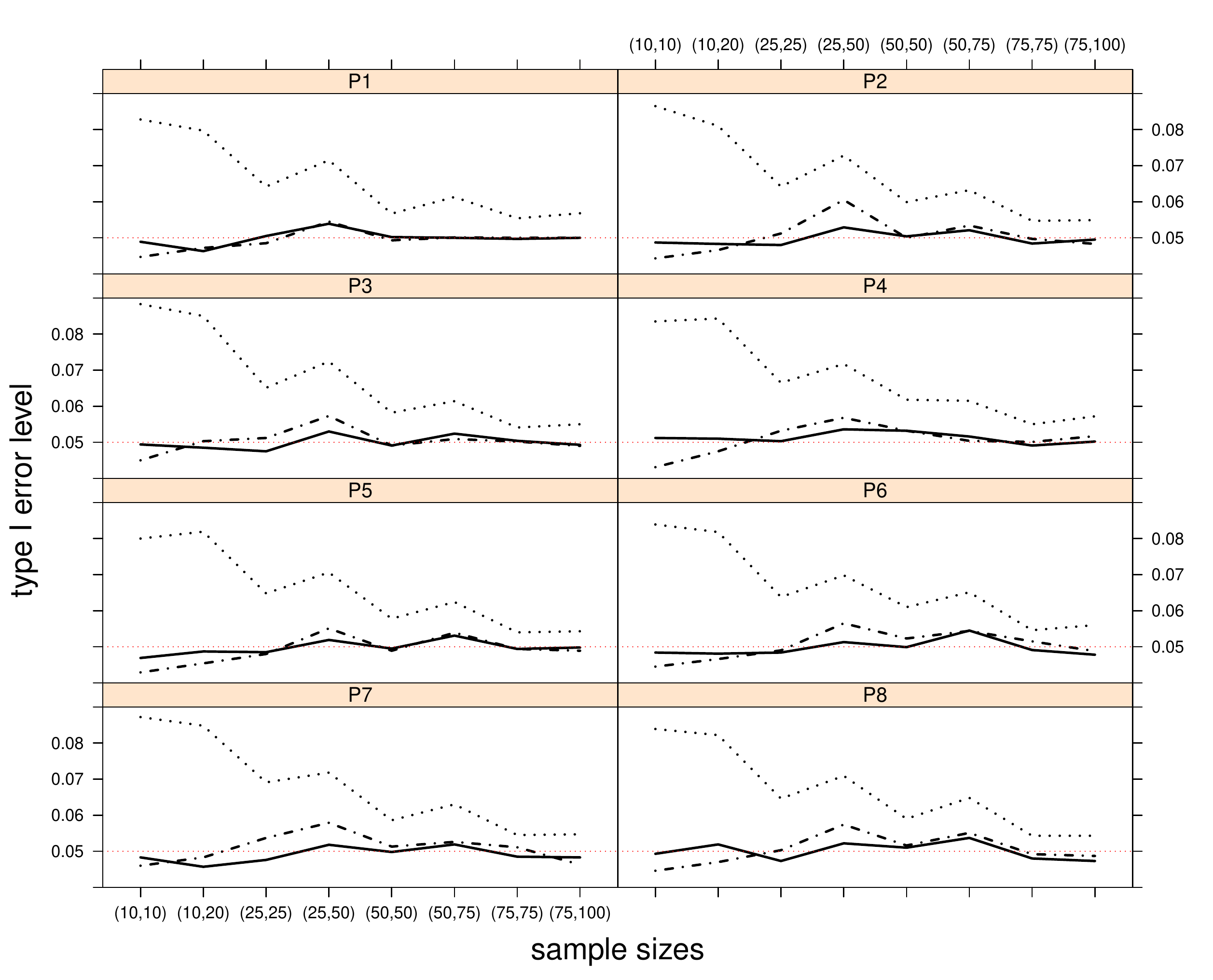}
\caption[kein Strich]{\label{fig:erg20_2}Type I error level ($\alpha=5\%$)  simulation results (y-axis) for threshold $\bs\tau_2$ and 20 items of the permutation test $\psi_n$
  (\rule[0.5ex]{0.7cm}{1pt}), the asymptotic test $\varphi_n$ (\hdashrule[0.5ex]{0.7cm}{1pt}{1pt}) and the bootstrap test $\psi_n^\star$ (\hdashrule[0.5ex]{0.7cm}{1pt}{1pt 1pt 5pt 1pt})
  for different sample sizes (x-axis).}
\end{figure}

As significance level $\alpha=0.05$ was chosen. The different statistical tests are compared by means of their type I error level.
The results for all simulations regarding five items are summarized in Figures \ref{fig:erg5_1} (threshold $\bs{\tau}_1$) and \ref{fig:erg5_2} (threshold $\bs{\tau}_2$). Figures \ref{fig:erg20_1} (threshold $\bs{\tau}_1$) and \ref{fig:erg20_2} (threshold $\bs{\tau}_2$) contain information about the 
simulations regarding twenty items.

The worst performing test in terms of the attained type I errors was the asymptotic test $\varphi_n$. Uniformly under all conditions of the 
simulation study. 
In contrast to the asymptotic test, the permutation and bootstrap tests performed reasonably well.
The statistical test regarding the permutation approach even yields better results than the statistical test based on the parametric bootstrap. Big differences between the
correlation matrices $\bs P_1, \ldots, \bs P_8$ are not observable. For small sample sizes, the permutation test seems to be the best choice for this testing procedure. The asymptotic test yields to liberal results, whereas the bootstrap method is slightly 
conservative. Differences between the number of items and the two different vectors of thresholds are not obvious. 
In light of these findings, the permutation test is the recommended procedure for testing or comparing two coefficient alpha for groups with an equal number of items. Whereas, the bootstrap test has its advantages in dealing with an unequal number
of items between groups and small to moderate sample sizes.
In particular, compared to an asymptotic approach, the permutation and bootstrap tests may be very useful, or indispensable in fact, 
for samples sizes in the range of $10$ to $100$ subjects. For larger sample sizes ($n_i>100$) the asymptotic testing provides a good control of the type-$I$-error rate and is thus recommended in these situations due to computational
efficiency. Additional simulations for continuous outcome given in the supplement, however, show that the ADF procedure may even have problems controlling the type I error rate for large sample sizes if data are rather skewed or 
the model assumptions are not fulfilled. Here the permutation procedure is even more advantageous.

\section{Application to empirical data}\label{sec:app}

This section is based on a work published by Maydeu-Olivares et al. (2010). The data example regarded in the latter publication is slightly extended in our work. 
The conducted negative problem orientation (NPO) questionnaire is one of five subscales of the Social Problem-Solving Inventory (SPSI-R, see D’Zurilla et al., 2002). A problem-solving ability was detected to be a process variable in several 
psychological disorders. Currently, two types of the NPO-questionnaire are available; a long- and a short-format test. The long form consists of ten items, whereas the short questionnaire only comprises 5 items. Each item is to be answered using a 
five-point response scale. Maydeu-Olivares et al. (2010) use two random samples from the U.S. population including 100 male and 100 female participants. The raw data is provided as supplementary material of the latter work.

Three examples were examined. The first one is the comparison of the reliability of the NPO questionnaire regarding men and women independently. And in a second example, the reliability of the long- and short-format is compared. 
In the third example the testing procedure for paired data is used. For a sample of overall 138 male and female participants the short-format test was repeated two times. Overall 10,000 permutation 
and 10,000 bootstrap samples were computed to calculate the critical value. The p-values of the permutation and bootstrap tests were calculated as described at the end of Section \ref{sec:perm}.

Regarding the subgroup analysis, coefficient alpha of the male sample is 0.837 and for the female sample it is 0.882. Consequently, the difference between these two coefficients is -0.045 (male -- female). The calculated test statistic $T_n$ 
has a value of -1.517 and yields to a p-value of 0.1291 regarding the asymptotic test. This result confirms the result in Maydeu-Olivares et al. (2010). The permutation test yields to a p-value of 0.1304. The permutation confidence interval for the 
corresponding test is given by $\left[-0.1034; 0.0129\right]$.

For the second example (long format vs. short format) only the bootstrap test is regarded since the permutation test is not valid for different test lengths. To ensure the independence of the observations of the groups, the long format data is based 
on the male population
and the data for the short format is based on all female participants. The calculated coefficient alpha for the long format test has a value of 0.837, whereas the value 
of the short format test is 0.776. Therefore, the difference between these values is 0.061 (male -- female). The corresponding test statistic $T_n$ has a value of 1.412 and the resulting p-value is 0.1580. The bootstrap test shows a 
one-sided p-value of 0.1978 and a two-sided one of 0.3956. These two values are in line with the results reported by Maydeu-Olivares et al. (2010).
Additionally to the p-values given above, a confidence interval
for the difference is calculated: $\left[-0.0502; 0.1713\right]$.

In the final example (repeated measures), the alpha coefficient for the first survey is $0.753$ and for the second  a value of $0.840$ is given. 
The test statistic $T_n$ from Section \ref{sec:boot} has a value of 
$-2.195$ which yields to a p-value of $0.0281$ for the asymptotic test. For the parametric bootstrap approach outlined in Section \ref{sec:boot} a p-value of $0.0384$ was calculated. Inverting this test, the corresponding bootstrap confidence 
interval for the difference is given by $\left[-0.1488;-0.0252\right]$.

\section{Conclusion}\label{sec:dis}

Cronbach (1951) coefficient alpha is \textit{the} reliability measure used by substantial or applied researchers. Whereas the single coefficient case has been extensively treated in literature, the comparison of two (or more)  coefficient alpha from 
a methodological statistical viewpoint has remained rather unstudied. 
In this paper, we have developed and discussed resampling-based and asymptotic tests for the two-sample testing problem:
$H_0 : \alpha_{C,1} = \alpha_{C,2}$ vs. $H_{11} : \alpha_{C,1} > \alpha_{C,2}$ (one-sided) or $H_{12} : \alpha_{C,1} \neq \alpha_{C,2}$ (two-sided) and outlined their extension to more general one-way designs.
We have investigated modified variants of permutation and bootstrap tests under the two-group general and most current \textit{asymptotically distribution-free} (ADF) framework and also discussed the more simple multivariate normal model. 
We have reported the results of a simulation study, which compared the resampling and asymptotic tests with regard to the type I error as the evaluation criterion. 
In particular, we have seen that the resampling techniques, and here especially the permutation tests, do improve the finite sample properties of the ADF asymptotic test, especially in small to moderate groups.

Although Cronbach's internal consistency coefficient alpha has been criticized on various grounds, this measure is still very popular and the most widely used reliability estimate.
In comparing groups, the standardization aspect of alpha (or other measures of reliability), deserves a
cautionary note. The comparison is based on the proportion (or ratio) of true-score variance relative to the observed total
variance, and it is in this relative sense that groups are compared. In particular, it may be possible, in the absolute or
unstandardized case, that the true-score variance (numerator) is equal in two groups, but the observed total variances
(denominator) may differ. This is a caveat indeed, and scientists must be aware that the share of explained variance is
what matters and is of interest here. The current theory may nevertheless be applied for comparisons of 
unstandardized effects, as long as they are given as adequately smooth functions of the underlying covariances. 

Another well-known limitation of Cronbach's alpha is that it is generally only a lower bound on the true reliability. 
Thus, it is actually tested whether the two groups have same or different lower bounds. 
In particular, the two groups may very well have identical reliability,
albeit the alpha bounds may differ. Or, it could also be the case that the two groups may have unequal reliabilities,
where the alpha lower bound may not differ. This an additional caveat to be aware of, of any lower bound reliability estimate.

The resampling-based inference methods presented in this paper 
and exemplified with coefficient alpha can also be applied to those alternative coefficients as mentioned in Section 4. Nevertheless, future research into this issue is needed. 

In particular, the analysis of more adequate reliability measures; especially for ordinal data, will be part of 
future research.
Another interesting direction for future research is to study the performance of the permutation and bootstrap tests for multiple (more than two) groups in extensive simulations. 
In this case, multiple comparisons and type I error inflation adjustment procedures may be of relevance and elaborated, too.  
Important from a practical viewpoint, the resampling-based and asymptotic tests have to be investigated and compared in more real applications than we have done. For example, in school comparison didactical surveys, samples typically are of school class
size, of approximately 15 to 35 pupils. In such practical situations, testing and comparing two or more coefficient alpha, or reliability estimates, based on the permutation tests may yield better results, i.e., exact and empirically valid conclusions. 
Applied future work may systematically explore analyses and comparisons of the resampling-based and asymptotic inference methods in realistic contexts and empirical datasets.





\newpage

\appendix

\section{Mathematical Appendix}
\setcounter{equation}{9}
Let the notation and prerequirements be as in Sections \ref{sec:intro} and \ref{sec:mod}. Also see, e.g., Muirhead (2009) and van der Vaart (1998) for the following multivariate and asymptotic elaborations.

\subsection*{ADF Asymptotics}

Let $\vec()$ be the usual operator that writes the elements of a symmetric matrix on and below the diagonal into a column vector, see e.g.  Muirhead (2009).
Due to the assumption of finite eighth order moments (ADF) and $\Ex(\bX_1)={\bf 0}$ we can write the normalized Cronbach coefficient as
\begin{equation}\label{aslinSigma}
  \vec\left(\sqrt{n_1}(\widehat{\bSigma}_1 - \bSigma_1)\right) =
  \frac{1}{\sqrt{n_1}}\sum_{i=1}^{n_1} \vec\left(\bX_i\bX_i' - \Ex(\bX_i\bX_i')\right) + o_p(1),
\end{equation}
where $o_p(1)$ converges in probability to zero as $n_1\to\infty$. Thus, it follows from the multivariate central limit theorem that
$\vec(\sqrt{n_1}(\widehat{\bSigma}_1 - \bSigma_1))$ is asymptotically multivariate normal with mean $\bf 0$ and covariance
$\Cov(\vec(\bX_1\bX_1'))$. Since $\alpha_{C,1} = \alpha_{C}(\bSigma_1)$ is a differentiable function of $\bSigma_1$ (or $\vec(\bSigma_1)$ respectively) it follows as in
Maydeu-Olivares et al. (2007) that $\sqrt{n_1}\widehat{\alpha}_{C,1}$ is asymptotically
normal distributed with mean $\alpha_{C,1}$ and variance $\tilde{\sigma}_1^2$ which depends on moments of fourth order.
In particular, the limit variance is given by
\begin{equation*}
  \tilde{\sigma}_1^2 = \tilde{\sigma}_1^2(\bSigma_1) = {\bs{\delta}(\bSigma_1)}'\var(\vec(\bX_1)){\bs{\delta}(\bSigma_1)},
\end{equation*}
which can be obtained from the delta method, see Maydeu-Olivares et al. (2007) for details.
Here the vector $\bs{\delta}(\bSigma_1)$ is a function of $\bSigma_1$ and is given in Equation~(4) in  Maydeu-Olivares et al. (2007).
However, we even know more.
Note, that $\alpha_{C,1}$ (as a function from $\R^{q_1}$ to $\R$, $q_1=\frac{k_1(k_1+1)}{2}$) is differentiable at $\vec(\bSigma_1)$ with
total derivative, i.e. Jacobi matrix, $\alpha_{\bSigma_1}'$, see van der Vaart (1998) for its explicit formula.
Hence, it follows from the proof of the multivariate delta method (to be concrete: the multivariate Taylor theorem),
see e.g. Theorem~3.1. in van der Vaart (1998), that $\widehat{\alpha}_{C,1} = \alpha_{C}(\widehat{\bSigma}_1)$ is even asymptotically linear
in this case, i.e.
\begin{equation}\label{aslin}
  \sqrt{n_1}(\widehat{\alpha}_{C,1} - \alpha_{C,1}) = \frac{1}{\sqrt{n_1}} \sum_{i=1}^{n_1} f_{\bSigma_1}(\bX_i) + o_p(1)
\end{equation}
holds as $n_1\to\infty$ with
\[f_{\bSigma_1}(\bX_i) = \alpha_{\bSigma_1}'\cdot\vec\left(\bX_i\bX_i' - \Ex(\bX_i\bX_i')\right).\]
The latter fulfills $\Ex(f_{\bSigma_1}(\bX_i))=0$ and $\var(f_{\bSigma_1}(\bX_i))=\tilde{\sigma}_1^2$.\\

Since a similar representation holds for $\sqrt{n_2}(\widehat{\alpha}_{C,2}-\alpha_{C,2})$ (with different variance $\tilde{\sigma}_2^2 = \tilde{\sigma}_2^2(\bSigma_2)$) it follows that
the statistic $M_n$ is also asymptotically normal under $H_0: \alpha_{C,1}=\alpha_{C,2}$ with mean zero
and variance $\tilde{\sigma}^2  = (1-\kappa)\tilde{\sigma}_1^2 + \kappa \tilde{\sigma}_2^2$, i.e.
\begin{align*}
  M_n
  &= \sqrt{\frac{n_2}{N}}\sqrt{n_1}(\widehat{\alpha}_{C,1}-\alpha_{C,1}) - \sqrt{\frac{n_1}{N}}\sqrt{n_2}(\widehat{\alpha}_{C,2}-\alpha_{C,2})\\
  &\dist \mathcal{N}(0,(1-\kappa)\tilde{\sigma}_1^2 + \kappa \tilde{\sigma}_2^2)
\end{align*}
if $n_1/N\rightarrow\kappa\in(0,1)$. A consistent estimator for $\tilde{\sigma}^2$ is given by
\begin{align}\label{pooled variance nonormal}
\begin{split}
 \widetilde{\sigma}^2
  =& \frac{n_2}{N} \left(\frac{1}{n_1-1} \sum_{i=1}^{n_1}\left(\widehat{\bs{\delta}}'(\bs{S}_{i1} - \bs{S}_1) \right)^2 \right)\\
  &+ \frac{n_1}{N}\left(\frac{1}{n_2-1} \sum_{i=1}^{n_2}\left(\widehat{\bs{\delta}}'(\bs{S}_{i2} - \bs{S}_2) \right)^2 \right),
\end{split}
\end{align}
see Equation~(7) in  Maydeu-Olivares et al. (2007)  for a similar formula in the one-sample case. Here, $\bs{S}_k = \vec(\widehat{\bSigma}_k)$ for $k=1,2$ 
and $\bs{S}_{i1} = \vec\left[(\bX_i - \overline{\bX}^{(1)})(\bX_i - \overline{\bX}^{(1)})'\right]$ for $1\leq i \leq n_1$ and $\overline{\bX}^{(1)} = \frac{1}{n_1} \sum_{i=1}^{n_1}\bX_i$ and $\bs{S}_{i2}$ is defined similarly
with the random variables of the second sample. Altogether it follows from Slutzky's theorem that the proposed
studentized test statistic $T_n =  T_n(\mathbb{X}) = \frac{M_n}{\widetilde{\sigma}}$ in \eqref{teststatnonnormal} is asymptotically standard normal under the null hypothesis
$H_0$, i.e. $T_n\dist \mathcal{N}(0,1)$.

\subsection*{Parametric Bootstrap}

To show that the proposed parametric bootstrap test $\psi_n^\star = \mathbbm{1}\{T_n > c_n^{\star}(\alpha)\}$ is of asymptotic level $\alpha$
we have to prove that the critical value $c_n^{\star}(\alpha)$, i.e. the conditional $(1-\alpha)$-quantile of the parametric bootstrap procedure,
converges in probability to the $(1-\alpha)$-quantile $z_{1-\alpha}$ of a standard normal distribution, i.e.
\begin{equation*}
 c_n^{\star}(\alpha) \prob z_{1-\alpha}
\end{equation*}
as $N\rightarrow\infty$, see Lemma~1 in Janssen \& Pauls (2003).
By continuity of the limit distribution, this is fulfilled if
the conditional parametric bootstrap distribution function of the test statistic $T_n$ is asymptotically standard normal in probability due to
$T_n\dist \mathcal{N}(0,1)$ under $H_0$. By assumption we again have
\begin{equation}\label{aslinSigmaboot}
 \vec\left(\sqrt{n_1}(\widehat{\bSigma}_1^\star - \bSigma_1)\right) =
  \frac{1}{\sqrt{n_1}}\sum_{i=1}^{n_1} \vec\left(\bX_i^\star\bX_i^{\star'} - \Ex(\bX_i^\star\bX_i^{\star'})\right) + o_p(1).
\end{equation}
Different to above, however, the family of random variables $\bX_i^\star, i\leq n_1$ now forms an array of row-wise i.i.d. random variables given the observed data. 
Thus, we cannot work with the classical multivariate CLT but have to employ the multivariate version of Lindeberg's or Lyapunov's theorems conditioned on the data. 
Due to the existence of finite eighth order moments and the consistency of $\widehat{\bSigma}_1$ 
Lyapunov's condition is fulfilled and we can obtain that
$\vec(\sqrt{n_1}(\widehat{\bSigma}_1^\star - \bSigma_1))$ is, given the data, asymptotically multivariate normal with mean $\bf 0$ and covariance matrix
$\Cov(\vec({\bf Z_1 \bf Z_1}'))$ in probability, where ${\bf Z_1}\sim \mathcal{N}(\bf 0, \bSigma_1)$.
Given the data, we can now proceed as in the prove above, i.e. we first apply the delta-method, then combine the results for
the two independent bootstrap samples and finally show that the given variance estimator is also consistent for the bootstrap (which follows, e.g.
from the Tchebyscheff inequality) to show that
\begin{equation*}
 \sup_{x\in \R}|P(T^\star \leq x|\bX_1,\dots,\bX_n) - \Phi(x)|\prob 0
\end{equation*}
as $\frac{n_1}{N}\rightarrow \kappa \in (0,1)$  and the result follows. Here, $\Phi$ is the distribution function of $\mathcal{N}(0,1)$. Due to the duality between statistical tests and confidence intervals, this also shows the asymptotic correctness of the 
latter. Moreover, the same argumentation also shows the lacking proof of the validity of Padilla et al. (2012) one-sample confidence interval for Cronbach's $\alpha$ coefficient.

\subsection*{Permutation Distribution}

Now suppose that $k_1=k_2$. In order to prove that the permutation test is of asymptotic level $\alpha$ we again have to show convergence of
the corresponding critical value $c_n^{\pi}(\alpha)$, i.e. the conditional $(1-\alpha)$-quantile of the permutation distribution function, converge in probability to
the $(1-\alpha)$-quantile $z_{1-\alpha}$ of a standard normal distribution, i.e.
\begin{equation*}
 c_n^{\pi}(\alpha) \prob z_{1-\alpha}.
\end{equation*}

In order to prove this, we apply Theorem~2.2 in Chung \& Romano (2013) together with a conditional Slutzky-type argument.

As in the beginning it holds that the normalized Cronbach coefficients are asymptotically linear in both groups, i.e. \eqref{aslin} as well as
\begin{equation}\label{aslin2}
  \sqrt{n_2}(\widehat{\alpha}_{C,2} - \alpha_{C,2}) = \frac{1}{\sqrt{n_2}} \sum_{i=n_1+1}^{N} f_{\bSigma_2}(\bX_i) + o_p(1),
\end{equation}
holds, where again $o_p(1)$ stands for a random variable that converges in probability to $0$ as $n_2\rightarrow\infty$.
Since by assumption $\tilde{\sigma}_1^2\in(0,\infty)$ all ingredients for applying Theorem~2.2 in  Chung \& Romano (2013) are fulfilled 
and it follows by Slutzky that
$$\frac{1}{N!} \sum_{\pi} \mathbbm{1}\{T_n(\mathbb{X}^\pi)\leq x\}$$
converges in probability to $\Phi(x)$. Altogether this proves that $\psi_n$ is an asymptotically exact level $\alpha$
testing procedure in the general ADF model.

\subsection*{Derivations for other reliability measures}

In the following the derivatives of the different reliability measures $\lambda_\ell, \ \ell=1,\dots, 8$ summarised in Section 4 of the main manuscript are given. Let $\bs \sigma=\vec(\bSigma)=\vec((\sigma_{ij})_{i,j})$, where $\vec()$ is a 
function stacking the elements of a symmetric matrix on and below the diagonal into a vector. Let $\bs \delta_\ell=\bs \delta_\ell(\bSigma)=\lambda'_\ell=\frac{d \lambda_\ell}{d \bs \sigma}$ be the derivative of $\lambda_\ell, \; \ell=1,\ldots, 8$. Below 
the entries of $\bs\delta_\ell$ are given:
\begin{align*}
 \frac{\partial \lambda_1}{\partial \sigma_{ij}} &= \left\{\begin{array}{ll} \frac{\tr(\bSigma)-\kone'\bSigma\kone}{\left(\kone'\bSigma\kone\right)^2}, & i=j \\[1em]
         2\cdot \frac{\tr(\bSigma)}{\left(\kone'\bSigma\kone\right)^2}, & i \neq j\end{array}\right. \\
 \frac{\partial \lambda_2}{\partial \sigma_{ij}} &= \left\{\begin{array}{ll} \frac{\tr(\bSigma)-\kone'\bSigma\kone-\sqrt{\frac{k}{k-1}}C_2^{\nicefrac{1}{2}}}{\left(\kone'\bSigma\kone\right)^2}, & i=j \\[1em]
         2\cdot \frac{1-\sqrt{\frac{k}{(k-1)}}\left(C_2^{\nicefrac{1}{2}}-\sigma_{ij}C_2^{-\nicefrac{1}{2}}\right)+\tr(\bSigma)-\kone'\bSigma\kone}{\left(\kone'\bSigma\kone\right)^2}, & i \neq j\end{array}\right. \\
 \frac{\partial \lambda_3}{\partial \sigma_{ij}} &= \left\{\begin{array}{ll} (-2)\cdot \frac{\kone'\bSigma\kone-\kone'\bSigma_A\kone-\kone'\bSigma_B\kone}{\left(\kone'\bSigma\kone\right)^2}, & i=j \\[1em]
         (-4)\cdot \frac{\kone'\bSigma\kone-\kone'\bSigma_A\kone-\kone'\bSigma_B\kone}{\left(\kone'\bSigma\kone\right)^2}, & i \neq j\end{array}\right.\\
 \frac{\partial \lambda_4}{\partial \sigma_{ij}} & = \left\{\begin{array}{ll} \frac{\tr(\bSigma)-\kone'\bSigma\kone+2\bar{C}_2^{\nicefrac{1}{2}}}{\left(\kone'\bSigma\kone\right)^2}, & i=j \\[1em]
         2\cdot \frac{\tr(\bSigma)-2\sigma_{ij}\left(\kone'\bSigma\kone\right)\bar{C}_2^{-\nicefrac{1}{2}}+2\bar{C}_2^{\nicefrac{1}{2}}}{\left(\kone'\bSigma\kone\right)^2}, & i \neq j\end{array}\right. \\
 \frac{\partial \lambda_{5}}{\partial \sigma_{ij}} &= \left\{\begin{array}{ll} \frac{\tr(\bSigma)-\kone'\bSigma\kone+\frac{2k}{k-1}\bar{C}_2^{\nicefrac{1}{2}}}{\left(\kone'\bSigma\kone\right)^2}, & i=j \\[1em]
        \frac{2\tr(\bSigma)-\frac{4k}{k-1}\left(\sigma_{ij}\left(\kone'\bSigma\kone\right)\bar{C}_2^{-\nicefrac{1}{2}}+\bar{C}_2^{\nicefrac{1}{2}}\right)}{\left(\kone'\bSigma\kone\right)^2}, & i \neq j\end{array}\right. \\
 \frac{\partial \lambda_{6}}{\partial \sigma_{ij}} &= \left\{\begin{array}{ll} \frac{\sum_{t=1}^k e_t^2}{\left(\kone'\bSigma\kone\right)^2}, & i=j \\[1em]
        2\cdot\frac{\sum_{t=1}^k e_t^2}{\left(\kone'\bSigma\kone\right)^2}, & i \neq j\end{array}\right.
\end{align*}

To handle the coefficients $\lambda_7$ and $\lambda_8$, we assume that $\bs c$ and $\bs A$ are differentiable in $\bSigma$ and additionally, we define the two differentiable functions $g(\bs c)=\kone'\bs{cc}'\kone$ and $h(\bs{A})=\kone'\bs{AA}'\kone$. Using the chain rule, the derivatives of $\lambda_7$ and $\lambda_8$ are given as follows:

\begin{align*}
 \frac{\partial \lambda_{7}}{\partial \sigma_{ij}} &= \left\{\begin{array}{ll} \frac{\left(\sum_{s=1}^k\frac{\partial g}{\partial  c_s}\frac{\partial c_s}{\partial \sigma_{ij}}+\sum_{s,t=1}^k \frac{\partial h}{\partial A_{st}}\frac{\partial A_{st}}{\partial \sigma_{ij}}\right) \left(\kone'\bSigma\kone\right)-\left(\kone'\bs{cc}'\kone+\kone'\bs{AA}'\kone\right) }{\left(\kone'\bSigma\kone\right)^2}, & i=j \\[1em]
         \frac{\left(\sum_{s=1}^k\frac{\partial g}{\partial c_s}\frac{\partial c_s}{\partial \sigma_{ij}}+\sum_{s,t=1}^k \frac{\partial h}{\partial A_{st}}\frac{\partial A_{st}}{\partial \sigma_{ij}}\right) \left(\kone'\bSigma\kone\right)-2\cdot\left(\kone'\bs{cc}'\kone+\kone'\bs{AA}'\kone\right) }{\left(\kone'\bSigma\kone\right)^2}, & i \neq j\end{array}\right.\\
  \frac{\partial \lambda_{8}}{\partial \sigma_{ij}} &= \left\{\begin{array}{ll} \frac{\sum_{s=1}^k\frac{\partial g}{\partial c_s}\frac{\partial c_s}{\partial \sigma_{ij}} \left(\kone'\bSigma\kone\right)-\left(\kone'\bs{cc}'\kone\right) }{\left(\kone'\bSigma\kone\right)^2}, & i=j \\[1em]
         \frac{\sum_{s=1}^k\frac{\partial g}{\partial c_s}\frac{\partial c_s}{\partial \sigma_{ij}} \left(\kone'\bSigma\kone\right)-2\cdot\left(\kone'\bs{cc}'\kone\right) }{\left(\kone'\bSigma\kone\right)^2}, & i \neq j\end{array}\right.
\end{align*}

Since $\vec\left(\sqrt{n}\left(\widehat{\bSigma}-\bSigma\right)\right)\stackrel{d}{\rightarrow} \mathcal{N}(\bs 0, \var(\vec(\bX_1 \bX_1')))$, where $\widehat{\bSigma}$ is the sample covariance matrix of independent and identically distributed random vectors 
$\bX_1, \ldots, \bX_n$ with $\bSigma=\Cov(\bX_1)$ and finite fourth moments, it thus, follows from the multivariate delta method that \[\sqrt{n}(\lambda_\ell(\widehat{\bSigma})-\lambda_\ell(\bSigma))\stackrel{d}{\rightarrow}\mathcal{N}(\bs 0, \bs\delta'_\ell\var(\vec(\bX_1 \bX_1'))\bs\delta_\ell)\]
for all choices of $\ell=1,\dots,8$. 
Due to the form of the derivatives given above the unknown variance can be consistently estimated; in case of $\lambda_7$ and $\lambda_8$ they depend on the specific forms of $\bs c$ and $\bs A$.

\section{More simulation results}
  Some simulation results for continuous data are presented. Two different scenarios are conducted: t-distributed and lognormally distributed data. In this section, we compare the ADF method to the permutation test presented in Section 3.1
  of the paper. The parametric bootstrap procedure has been left out for lucidity and since the permutation method performed slightly better. Moreover, recall that the permutation test is finitely exact under exchangeability. Again 10,000 simulation trails with 500 
  permutation samples were performed.

  To check the behavior of the procedures in case of deviations from the underlying moment assumption, we first deal with t-distributed data with four degrees of freedom. Note, that the assumption of finite eight order moment is clearly violated in this case. 
  The data are generated with the help of the \textsc{R} function
  \texttt{rmvt()} which is included in the \texttt{mvtnorm} package. Based on the simulation results of the main manuscript, the following results are based on two correlation matrices only. The reason is that there are matrices following the true-score 
  equivalent model and some do not. Another cause is the comparability of the results of the main simulation study. Thus, in the following we only consider correlation matrices $\bs{P_1}$ and $\bs{P_4}$ given in Section 5.1. 

\begin{figure}[!ht]
  \centering
  \subfigure[equal sample sizes]{\includegraphics[width=0.45\textwidth]{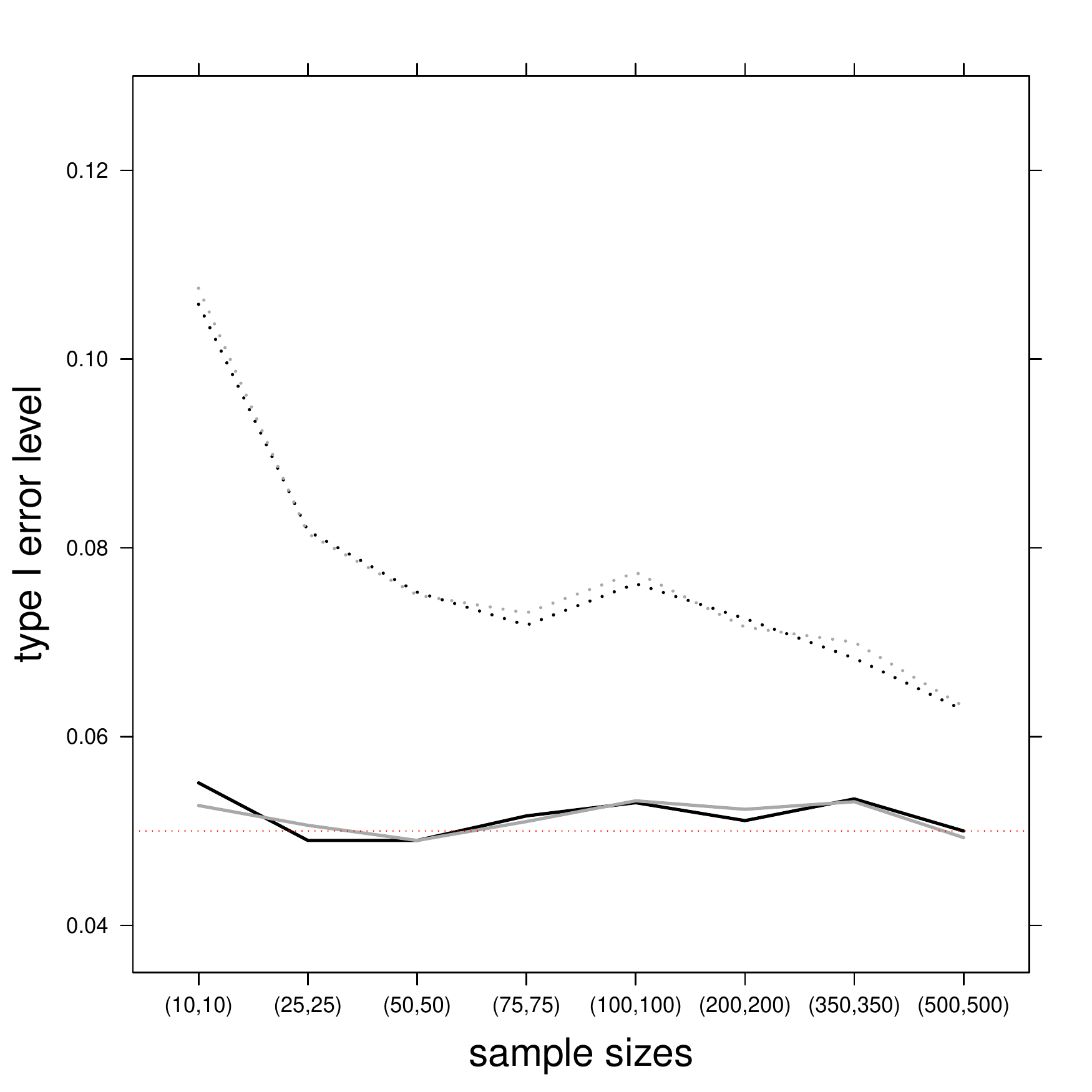}}\qquad
  \subfigure[unequal sample sizes]{\includegraphics[width=0.45\textwidth]{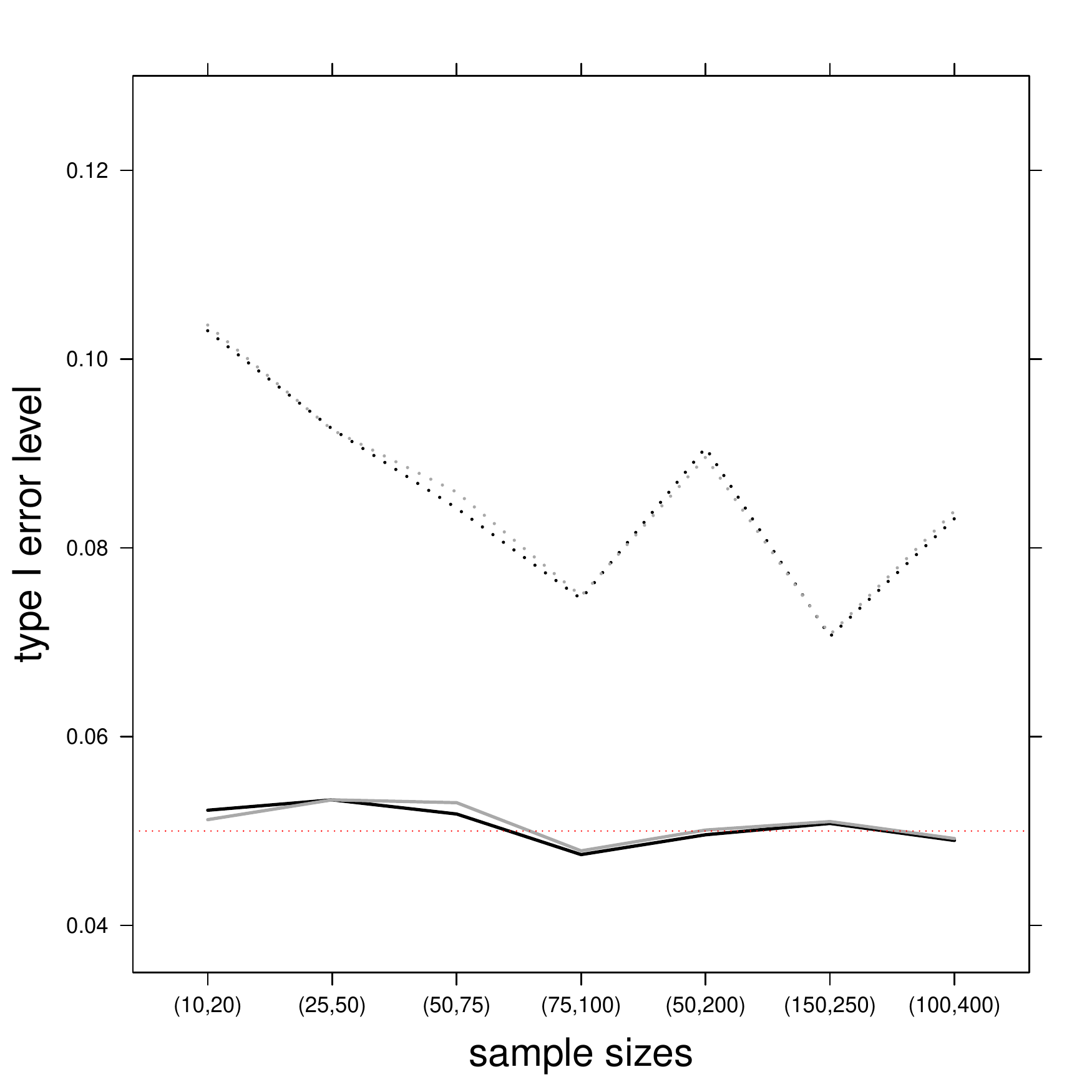}}\\
  \caption[text]{\label{fig:tdistr}Type I error level ($\alpha=5\%$)  simulation results (y-axis) for t-distributed data of the permutation test $\psi_n$
  (\rule[0.5ex]{0.7cm}{1pt}) and the asymptotic test $\varphi_n$ (\hdashrule[0.5ex]{0.7cm}{1pt}{1pt}) 
  for different sample sizes (x-axis) and two different correlation matrices $\bs{P_1}$ (black) and $\bs{P_4}$ (grey).}
\end{figure}

The results are summarized in Figure~\ref{fig:tdistr}, where the type I error levels of the permutation and the asymptotic test for two different correlation matrices are shown. In the left plot, same sample sizes in the different groups are considered, whereas the right plot summarizes the results
of unequal sample sizes. It is evident that the asymptotic test does not control the type I error rate satisfactorily in all cases. Even for very large balanced sample sizes ($n_i>350$) the type I errors are still around 7\% and even larger in extremely 
unbalanced cases or smaller sample sizes. 
In contrast, the novel permutation test controls the type I error rate fairly well in all situations and is always in the range of 4.7 and 5.3\%.

\begin{figure}[!ht]
  \centering
  \subfigure[equal sample sizes]{\includegraphics[width=0.45\textwidth]{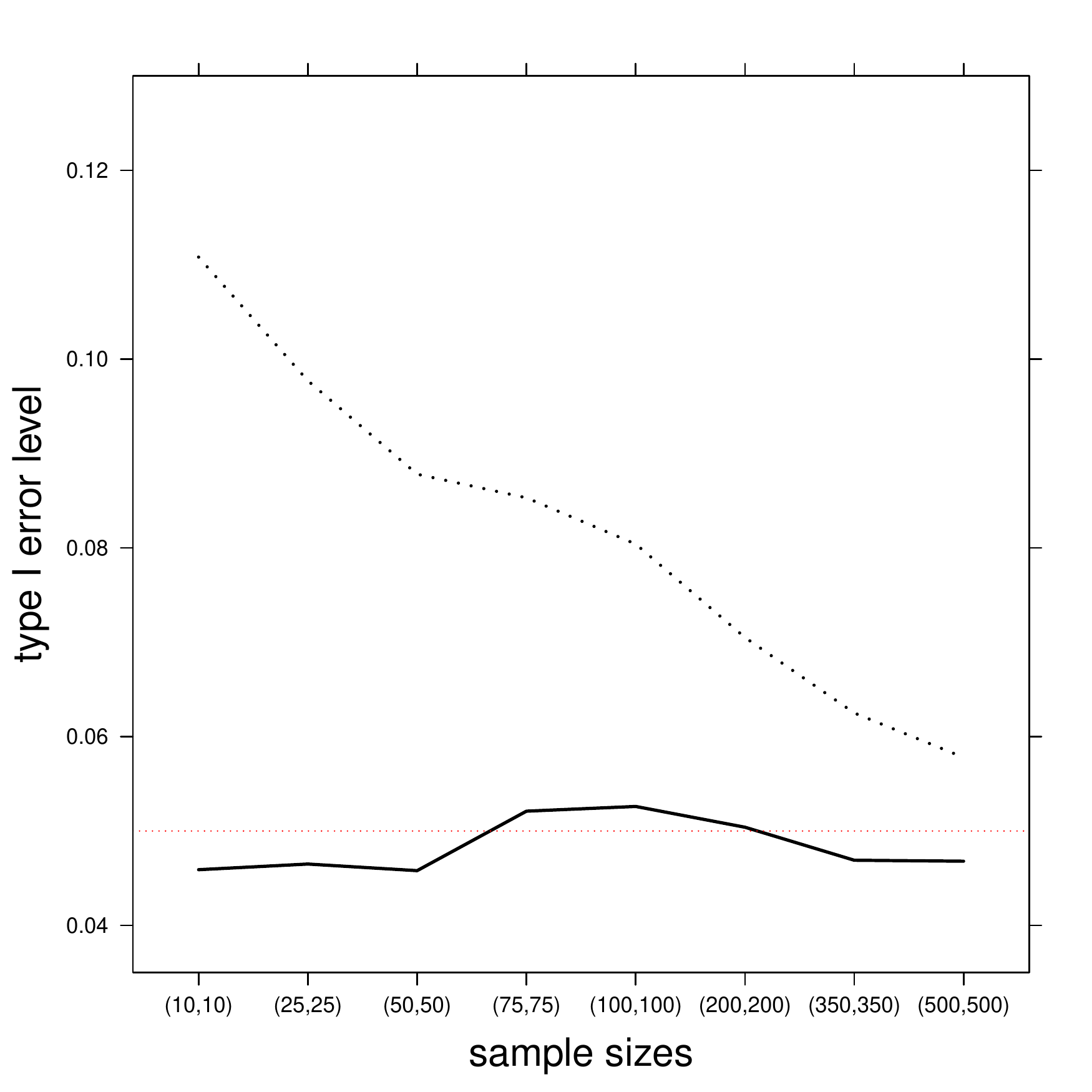}}\qquad
  \subfigure[unequal sample sizes]{\includegraphics[width=0.45\textwidth]{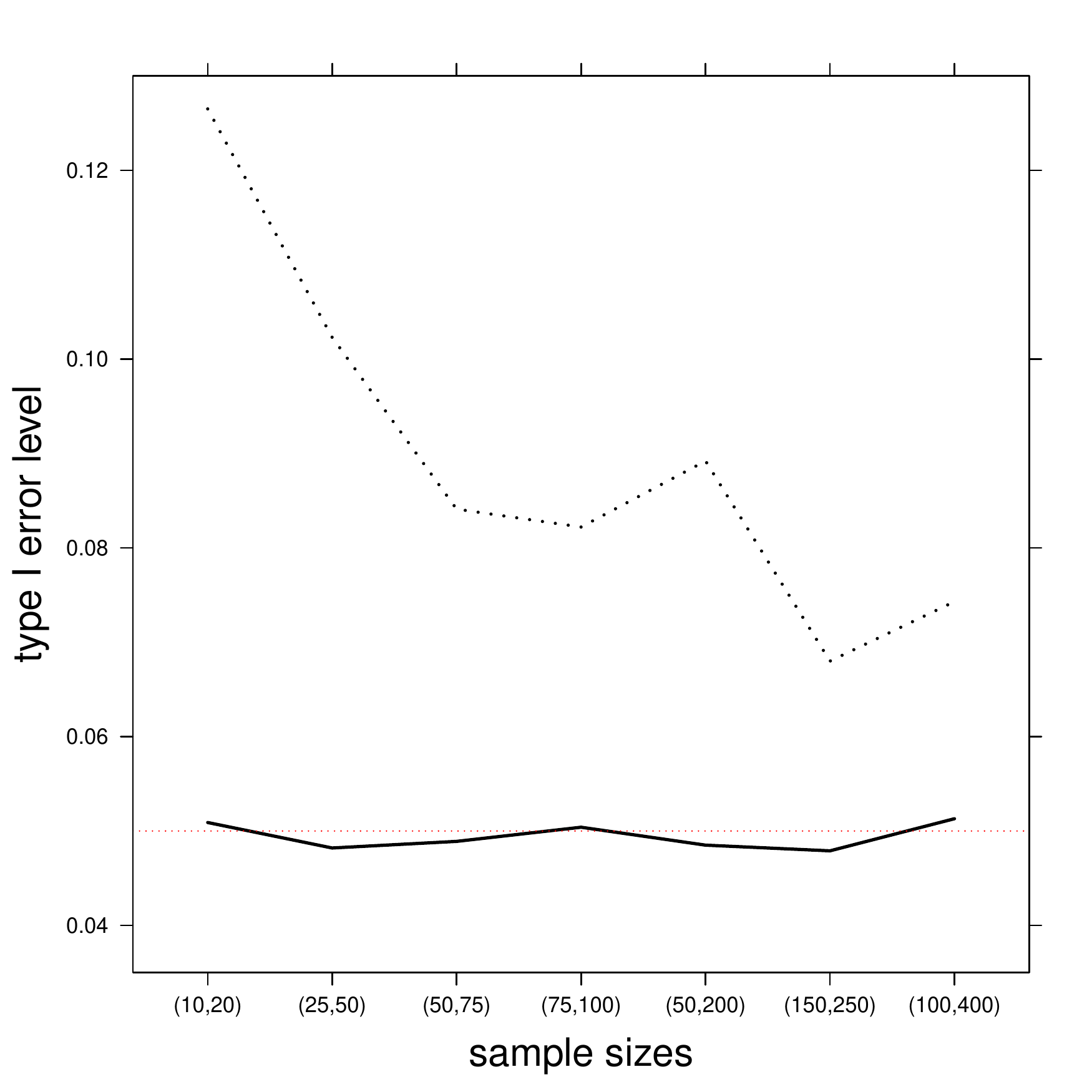}}\\
  \caption[text]{\label{fig:lognormdistr}Type I error level ($\alpha=5\%$)  simulation results (y-axis) for lognormal distributed data of the permutation test $\psi_n$
  (\rule[0.5ex]{0.7cm}{1pt}) and the asymptotic test $\varphi_n$ (\hdashrule[0.5ex]{0.7cm}{1pt}{1pt}) 
  for different sample sizes (x-axis).}
\end{figure}

Next, we deal with log-normally distributed data. The data are generated by a scale model with $k=5$ items $\bX_i= \bs{I}_k^{1/2}\bep_i, \; i=1,\ldots, N$, where $\bep_i=\frac{e_i-\Erw(e_i)}{\sqrt{\var(e_i)}}$ and $e_i\sim LN(0,1)$ are independent standardized log-normally distributed error terms. Figure~\ref{fig:lognormdistr} shows the results of the log-normal
distribution. Contrary to the situation with the t-distributed data before, observations simulated under this scenario fulfill the postulated moment assumption. However, the observations are rather similar. For smaller or strongly unbalanced sample sizes the true 
type I error is around 10\% (or even larger) and decrease with increasing $n_i$. However, even for larger sample sizes the type I error control is not very satisfactory. The asymptotic test exhibits some issues in controlling the type I error rate, whereas the permutation test works quite perfect.

\section{R code}
In the following, we present the \textsc{R} code of our new permutation and parametric bootstrap procedures. First, we present the different functions for the two resampling methods (\texttt{pval.perm()} and \texttt{pval.boot()}). In a third part, 
a function for calculating the test statistic (\texttt{tstat()}) and another
function which writes the elements of a symmetric matrix on and below the diagonal into a column vector (\texttt{vecs()}) are given. 
\subsection{R code of the permutation test}\label{sec:Rperm}

\begin{lstlisting}[language=R]
  pval.perm <- function(data, n1 = NULL, n2 = NULL, p = NULL, B = 1000){
    library(MASS)
    perm.results <- matrix(rep(0, (4 * B)), ncol = 4)
    n <- n1 +n2
    
    # original data estimates of alpha and T statistics
    orig.results <- tstat(data1, data2, n1, n2, p1, p2)
    
    # permuted data estimates of alpha and T statistics
    for (i in 1:B){    
      dat_temp <- data[sample(1:nrow(data)),]
      perm.results[i, ] <- tstat(dat_temp, n1, n2, p)
    }
  perm.p.values = perm.p.values_nonorm <- numeric(3)
  
  # permutation p-values
  perm.p.values[1] <- (sum(orig.results[1] <= perm.results[, 1]) / B)  # right-sided
  perm.p.values[2] <- (sum(orig.results[1] >= perm.results[, 1]) / B)  # left-sided
  perm.p.values[3] <- (2 * min(perm.p.values[1:2]))  # two-sided
  names(perm.p.values) <- c("right.sided", "left.sided", "two-sided")
  perm.p.values_nonorm[1] <- (sum(orig.results[2] <= perm.results[, 2]) / B)  # right-sided
  perm.p.values_nonorm[2] <- (sum(orig.results[2] >= perm.results[, 2]) / B)  # left-sided
  perm.p.values_nonorm[3] <- (2 * min(perm.p.values_nonorm[1:2]))  # two-sided
  names(perm.p.values_nonorm) <- c("right.sided", "left.sided", "two-sided")
  
  return(list(perm.p.values=perm.p.values, perm.p.values_nonorm=perm.p.values_nonorm, alpha1=orig.results[3], alpha2=orig.results[4]))
}
\end{lstlisting}

\subsection{R code of the parametric bootstrap test}\label{sec:Rboot}

\begin{lstlisting}[language=R]
pval.boot <- function(data, n1 = NULL, n2 = NULL, p = NULL, B = 1000){
  library(mvtnorm)
  boot.results <- matrix(rep(0, (2 * B)), ncol = 2)
  n <- n1+n2
  
  # original data estimates of alpha and T statistics
  orig.results <- tstat(data, n1, n2, p)
  
  # bootstraped data estimates of alpha and T statistics
  Sigma1 <- cov(data[1:n1, 1:p])
  Sigma2 <- cov(data[(n1 + 1):n, 1:p])
  for (i in 1:B){
    dat_temp <- rbind(mvrnorm(n1, rep(0, p), Sigma1), mvrnorm(n2, rep(0, p), Sigma2))
    boot.results[i, ] <- tstat(dat_temp, n1, n2, p)
  }  
  
   # bootstrap p-values
   boot.p.values[1] <- (sum(orig.results[3] <= boot.results[, 3]) / B)  # right-sided
   boot.p.values[2] <- (sum(orig.results[3] >= boot.results[, 3]) / B)  # left-sided
   boot.p.values[3] <- (2 * min(boot.p.values[1:2]))  # two-sided
   names(boot.p.values) <- c("right.sided", "left.sided", "two-sided")
   boot.p.values_nonorm[1] <- (sum(orig.results_nonorm[3] <= boot.results_nonorm[, 3]) / B)  # right-sided
   boot.p.values_nonorm[2] <- (sum(orig.results_nonorm[3] >= boot.results_nonorm[, 3]) / B)  # left-sided
   boot.p.values_nonorm[3] <- (2 * min(boot.p.values_nonorm[1:2]))  # two-sided
   names(boot.p.values_nonorm) <- c("right.sided", "left.sided", "two-sided")
  
  
  return(list(boot.p.values=boot.p.values, boot.p.values_nonorm=boot.p.values_nonorm, alpha1=orig.results[3], alpha2=orig.results[4]))
}
\end{lstlisting}

\subsection{R code of the test statistic and the vecs-function}\label{sec:Rother}

\begin{lstlisting}[language=R]
### function vecs
vecs <- function(data){
  upna <- data
  upna[upper.tri(data)] <- NA
  upna_vec <- as.vector(upna)[!is.na(as.vector(upna))]
  return(as.matrix(upna_vec))
}
\end{lstlisting}

\begin{lstlisting}[language=R]
### calculates the test statistics of both tests
tstat <- function(data, n1 = NULL, n2 = NULL, p = NULL){
  n <- (n1 + n2)
  Sigma1 <- cov(data[1:n1, 1:p])
  Sigma2 <- cov(data[(n1 + 1):n, 1:p])
  col.mean1 <- matrix(colMeans(data[1:n1, 1:p]), nrow = 1) 
  col.mean2 <- matrix(colMeans(data[(n1 + 1):n, 1:p]), nrow = 1) 
  trSigma1 <- sum(diag(Sigma1))
  trSigma2 <- sum(diag(Sigma2))
  sSigma1 <- sum(Sigma1)
  sSigma2 <- sum(Sigma2)
  
  # variances, separately
  sigma1q <- ((2 * p^2 * (sSigma1 * (sum(diag(Sigma1 %*% Sigma1)) + trSigma1^2) - 2 * trSigma1 * sum(Sigma1 %*% Sigma1))) / ((p - 1)^2 * sSigma1^3))
  sigma2q <- ((2 * p^2 * (sSigma2 * (sum(diag(Sigma2 %*% Sigma2)) + trSigma2^2) - 2 * trSigma2 * sum(Sigma2 %*% Sigma2))) / ((p - 1)^2 * sSigma2^3))
                          
  # Welch-type variance, pooled 
  sigma <- sqrt((n2 / n) * sigma1q + (n1 / n) * sigma2q)
  
  # variances nonorm, separately
  helpdelta1 <- 2*p/(p-1)*(trSigma1/(sSigma1)^2)
  helpdeltatr1 <- -p/(p-1)*((sSigma1-trSigma1)/(sSigma1)^2)
  delta_1 <- matrix(rep(helpdelta1, p^2), nrow = p)
  diag(delta_1) <- helpdeltatr1
  
  helpdelta2 <- 2*p/(p-1)*(trSigma2/(sSigma2)^2)
  helpdeltatr2 <- -p/(p-1)*((sSigma2-trSigma2)/(sSigma2)^2)
  delta_2 <- matrix(rep(helpdelta2, p^2), nrow = p)
  diag(delta_2) <- helpdeltatr2
  
  sigma1q.non <- 0
  wcv <- 0
  v <-0
  tmp <- 0
  for (i in 1:n1){
    v <- (as.matrix(data[i,1:p, drop = FALSE]) - col.mean1)
    wcv <- (t(vecs(delta_1))%*%(vecs((t(v) %*%v))-vecs(Sigma1)))^2
    sigma1q.non <- (sigma1q.non + wcv)
  }
  
  sigma2q.non <- 0
  wcv <- 0
  v <-0
  tmp <- 0
  for (i in 1:n2){
    v <- (as.matrix(data[i+n1,1:p, drop = FALSE]) - col.mean2)
    wcv <- (t(vecs(delta_2))%*%(vecs((t(v) %*%v))-vecs(Sigma2)))^2
    sigma2q.non <- (sigma2q.non + wcv)
  }
  
  # variance, pooled 
  sigma.non <- sqrt(n2/n*(1/(n1-1)*sigma1q.non)+n1/n*(1/(n2-1)*sigma2q.non))
  
  # Cronbach alpha
  alpha1 <- (p / (p - 1) * (1 - trSigma1 / sSigma1))
  alpha2 <- (p / (p - 1) * (1 - trSigma2 / sSigma2))
  
  # test statistic
  Mn <- (sqrt((n1 * n2) / n) * (alpha1 - alpha2))
  tval <- (Mn / sigma)
  tval.nonorm <- (Mn / sigma.non)
  
  return(c(TSTAT = tval, TSTAT_NONORM = tval.nonorm))
}
\end{lstlisting}

\vspace{\fill}


\end{document}